\documentclass[a4paper,leqno,10pt]{amsart}

\usepackage{epsf,graphicx,epsfig}
\usepackage{amscd}
\usepackage{amsmath,latexsym,amssymb,amsthm}
\usepackage[nospace,noadjust]{cite}
\usepackage{textcomp}
\usepackage{setspace,cite}
\usepackage{lscape,fancyhdr,fancybox}
\usepackage{stmaryrd}
\usepackage[all,cmtip]{xy}
\RequirePackage{tikz-cd}
\RequirePackage{amssymb}
\usetikzlibrary{calc}
\usetikzlibrary{decorations.pathmorphing}

\tikzset{curve/.style={settings={#1},to path={(\tikztostart)
	.. controls ($(\tikztostart)!\pv{pos}!(\tikztotarget)!\pv{height}!270:(\tikztotarget)$)
	and ($(\tikztostart)!1-\pv{pos}!(\tikztotarget)!\pv{height}!270:(\tikztotarget)$)
	.. (\tikztotarget)\tikztonodes}},
	settings/.code={\tikzset{quiver/.cd,#1}
	\def\pv##1{\pgfkeysvalueof{/tikz/quiver/##1}}},
	quiver/.cd,pos/.initial=0.35,height/.initial=0}

\tikzset{tail reversed/.code={\pgfsetarrowsstart{tikzcd to}}}
\tikzset{2tail/.code={\pgfsetarrowsstart{Implies[reversed]}}}
\tikzset{2tail reversed/.code={\pgfsetarrowsstart{Implies}}}
\tikzset{no body/.style={/tikz/dash pattern=on 0 off 1mm}}

\usepackage{graphicx}

\usepackage{color}
\usepackage{url}
\usepackage{enumerate}
\usepackage[mathscr]{euscript}

\newcounter{cnt1}
\newcounter{cnt2}
\newcounter{cnt3}
\newcommand{\blr}{\begin{list}{$($\roman{cnt1}$)$} {\usecounter{cnt1}
		\setlength{\topsep}{0pt} \setlength{\itemsep}{0pt}}}
	    \newcommand{\bla}{\begin{list}{$($\alph{cnt2}$)$} {\usecounter{cnt2}
		\setlength{\topsep}{0pt} \setlength{\itemsep}{0pt}}}
		\newcommand{\bln}{\begin{list}{$($\arabic{cnt3}$)$} {\usecounter{cnt3}
		\setlength{\topsep}{0pt} \setlength{\itemsep}{0pt}}}
		\newcommand{\el}{\end{list}}
		
		\newtheorem{Thm}{Theorem}[section]
		
		\newtheorem{Prop}[Thm]{Proposition}
		\newtheorem{Def}[Thm]{Definition}
		\newtheorem{Exm}[Thm]{Example}
		\newtheorem{Rem}[Thm]{Remark}
		\newtheorem{Cor}[Thm]{Corollary}

		\title{}
		\author{}
		\date{}
		\usepackage{amssymb}

		\usepackage{hyperref}
		\hypersetup{
			colorlinks,
			citecolor=blue,
			filecolor=black,
			linkcolor=blue,
			urlcolor=black
		}
	
\begin{document}
			
\title{Lie algebroids, quantum Poisson algebroids, and Lie algebroid connections}
\author{Satyendra Kumar Mishra and Abhishek Sarkar}
\begin{abstract}
In this paper, we consider Lie algebroids over commutative ringed spaces. Lie algebroids over ringed spaces unify the existing notion of Lie algebroids over smooth manifolds, complex manifolds, analytic spaces, algebraic varieties, and schemes. We show that the universal enveloping algebroid of a Lie algebroid possesses a natural filtration that yields a structure of a sheaf of quantum Poisson algebras. We establish a bijective correspondence between sheaves of quantum Poisson algebras and Lie algebroids. We show that this correspondence leads to an adjunction between the two categories. We discuss this bijective correspondence in particular cases of Lie algebroids over ringed spaces and highlight the subsequent results. To characterize non-flat Lie algebroid connections, we construct a sheaf of twisted universal enveloping algebras for a Lie algebroid using Lie algebroid (hyper) cohomology. We show that our construction yields some of the existing constructions for Lie-Rinehart algebras and holomorphic Lie algebroids. As another application, we study the deformation groupoid of a Lie algebroid using the second hypercohomology of the Lie algebroid. 
\end{abstract}
\footnote{AMS Mathematics Subject Classification $(2020)$: $17$B$35$, $32$L$10$, $53$D$17$, $14$F$10$}
\keywords{ Lie algebroids, ringed spaces, Lie-Rinehart algebras,  universal enveloping algebra, quantum Poisson algebras, Lie algebroid connections, Lie algebroid hypercohomology}
\maketitle
\section{Introduction}

\noindent The concept of Lie algebroids over smooth manifolds plays a fundamental role in differential geometry and mathematical physics, as they generalize infinitesimal symmetries of spaces \cite{MM, RF}. These symmetries are closely linked to the corresponding global symmetries described by Lie groupoids \cite{KM}. 
An algebraic counterpart of Lie algebroids, known as Lie-Rinehart algebras, provides a framework to study more general situations \cite{JH, MK, LM, UB}. In complex and algebraic geometry, Lie algebroids have been studied over analytic spaces \cite{BP} and over algebraic varieties (or schemes) \cite{MK, CV, UB}, respectively. 

\medskip

\noindent In this article, we study Lie algebroids over commutative ringed spaces, see Definition \ref{Category LA} for more details. Several recent works, including those of Mikhail Kapranov~\cite{MK}, Michel van den Bergh et al.~\cite{CV,CRV, DRV}, and other authors ~\cite{SC, LG-PX, JV, AA, BP, BRT}
considered sheaf-theoretic analogs of Lie algebroids to discuss smooth, complex analytic, and algebraic cases uniformly. The sheaf-theoretic approach offers a unified framework for studying global calculus on both smooth and singular geometric spaces (see \cite{SR, BP, TS}). This perspective facilitates the analysis of diverse geometric structures, including Poisson analytic spaces \cite{BP}, singular foliations, generalized involutive distributions \cite{RF, BP, LG}, and free Lie algebroids \cite{MK}.

\medskip
            
\noindent The universal enveloping algebra of a Lie-Rinehart algebra plays a fundamental role in understanding the associated representation theory, deformation theory and homological algebra \cite{GR, MM, JH, TS}. The universal enveloping algebra of a Lie-Rinehart algebra has a natural filtration, which makes it a `quantum Poisson algebra'. In particular, consider the Lie-Rinehart algebra consists of derivations of the $\mathbb{R}$-algebra of smooth functions $C^{\infty}(M)$ on a manifold $M$ (see \cite{MM, JH}). In this case, the notion of quantum Poisson algebras naturally arises when studying its universal enveloping algebra, namely the algebra of differential operators $\mathcal{D} (M)$ (see \cite{GP, TS}). In \cite{GP}, Grabowski and Poincin studied automorphisms on $\mathcal{D} (M)$, with respect to the standard quantum Poisson algebra structure. In recent years, the universal enveloping algebra construction has been extended to the setting of Lie algebroids over ringed spaces (see \cite{CV, DRV, UB}). Let $\mathcal{L}$ be a Lie algebroid over a (commutative) ringed space $(X,\mathcal{O})$ or, equivalently, $(\mathcal{O}, \mathcal{L})$ be a Lie algebroid over a topological space $X$. If we take the associated presheaf of universal enveloping algebras and sheafifying it, we get the universal enveloping algebroid $\mathscr{U}(\mathcal{O},\mathcal{L})$  of the Lie algebroid $\mathcal{L}$.  Following the terminology of `universal enveloping algebroid', we call a sheaf of quantum Poisson algebras a `quantum Poisson algebroid'. 

\medskip

\noindent The first goal of this paper is to derive an adjunction between the category of Lie algebroids (over $X$) and the category of quantum Poisson algebroids (over $X$). To achieve this goal, we first show that the universal enveloping algebroid $\mathscr{U}(\mathcal{O},\mathcal{L})$ has a natural quantum Poisson algebroid structure. Thus, for any Lie algebroid $(\mathcal{O}, \mathcal{L})$ (over $X$), we can associate a sheaf of quantum Poisson algebras, namely the sheaf $\mathscr{U}(\mathcal{O},\mathcal{L})$. Conversely, we associate a Lie algebroid structure with a given sheaf of quantum Poisson algebras over $X$ (Theorem \ref{quant-Lie}). This correspondence yields an adjunction between the category of Lie algebroids and the category of quantum Poisson algebroids (Theorem \ref{LA-QP sheaf}). We also discuss this correspondence in the particular cases of Lie-Rinehart algebras, smooth and holomorphic Lie algebroids.  

\medskip

\noindent The second goal of this paper is to study twisted universal enveloping algebroids and connections of a Lie algebroid over a ringed space.  In \cite{Sridh}, Sridharan introduced the notion of twisted universal enveloping algebras of Lie algebras to study representation theory within the framework of filtered algebras. In a related development, Simpson introduced the concept of almost polynomial filtered algebras in \cite{Simpson}, aimed at constructing moduli spaces of representations of the fundamental group of smooth projective varieties. Building on these ideas, Tortella extended the framework to the holomorphic setting in \cite{PT}, studying moduli spaces of semistable flat $\mathcal{L}$-connections for a holomorphic Lie algebroid $\mathcal{L}$ using almost polynomial filtered algebras satisfying Simpson’s axioms. Later, in \cite{HM}, Maakestad further extended the theory to the setting of Lie-Rinehart algebras, where the notion of almost commutative rings was developed to study moduli spaces of sheaves of rings of differential operators over affine schemes.

\medskip

\noindent  To address the second goal, we consider the associated cochain complex, known as the Chevalley-Eilenberg-de Rham complex \cite{UB, BRT, BP, PT} and construct the twisted universal enveloping algebroid $\mathscr{U}(\mathcal{O},\mathcal{L}, \omega)$ of a Lie algebroid $(\mathcal{O}, \mathcal{L})$ together with a 2-cocycle $\omega \in \mathcal{Z}^2(\mathcal{L}, \mathcal{O})$ of the cochain complex. This object is equipped with the structure of a standard sheaf of almost commutative filtered algebras. Equivalently, we show that this construction defines a quantum Poisson algebroid. In a particular case, we show that our construction gives the twisted universal enveloping algebras of a Lie-Rinehart algebra studied by Maakestad in \cite{HM}. Moreover, in the holomorphic context, we obtain Tortella's construction as in \cite{PT}. We show that for an  $\mathcal{O}$-module $\mathcal{E}$, an $\mathscr{U}(\mathcal{O},\mathcal{L}, \omega)$-module structure on $\mathcal{E}$ is equivalent to an $\mathcal{L}$-connection on $\mathcal{E}$ with curvature type $\omega$, where $(\mathcal{O}, \mathcal{L})$ is a locally free Lie algebroid  (Theorem \ref{connection with prescribed curvature}).
 We also establish a version of the Poincaré-Birkhoff-Witt (PBW) theorem for $\mathscr{U}(\mathcal{O}, \mathcal{L},\omega)$ (Theorem \ref{twisted PBW theorem}). 
 Finally, we relate the cohomology group \( \mathbb{H}^2(\mathcal{L}, \mathcal{O}) \) to the deformation of sheaves of almost commutative filtered algebras from the classical to the quantum setting.
Subsequently, we study the deformation groupoid $\mathcal{A}(\mathscr{S}_{\mathcal{O}} \mathcal{L})$ of a locally free Lie algebroid $\mathcal{L}$, where $\mathscr{S}_{\mathcal{O}} \mathcal{L}$ is the sheaf of symmetric algebras of $\mathcal{L}$ over $\mathcal{O}$ (Theorem \ref{deformation groupoid}).

\medskip
 \noindent Beyond their role in representation theory, quantum Poisson algebroids and sheaves of almost commutative filtered algebras also emerge as fundamental structures in noncommutative geometry and mathematical physics. In particular, they play a central role in Poisson geometry, which provides a conceptual bridge to noncommutative geometry via deformation quantization and Hamiltonian mechanics (see \cite{TS}). Moreover, they are naturally encountered in the study of quantum groups, where Poisson modules over Poisson algebras are central (see \cite{HB}). We refer to \cite{VGinz2, LM} for further developments in the theory of $\mathcal{D}$-modules and related areas. In this work, we further explore some of these aspects in the contexts of Poisson geometry, logarithmic foliations, meromorphic connections, (logarithmic) de Rham cohomology, and Atiyah algebroids.

\medskip

 \noindent In Section \ref{Sec 2}, we recall standard notions related to Lie algebroids in the algebro-geometric setting, including universal enveloping algebroids, Atiyah algebroids, and the hypercohomology of Lie algebroids. In Section \ref{Sec 3}, we introduce the concept of quantum Poisson algebroids and explore their correspondence with Lie algebroids. Along the way, we discuss key examples such as universal enveloping algebroids and the sheaf of differential operators on an $\mathcal{O}$-module. We also examine special cases including Lie-Rinehart algebras, smooth and holomorphic Lie algebroids, and singular foliations. In Section \ref{Sec 4}, we study twisted universal enveloping algebroids associated with Lie algebroid connections, which naturally carry a quantum Poisson algebroid structure. We show that this construction unifies various existing constructions found in the literature. 
 
 \section{Preliminaries}  \label{Sec 2}
In this section, we recall the basic definitions and examples of Lie algebroids over ringed spaces and its universal enveloping algebroids. We also recall the relationship between these two notions. Then we recall some properties of the Chevalley-Eilenberg-de Rham hypercohomology of such a Lie algebroid.

\subsubsection{Notations and conventions.}There is a canonical procedure to construct a sheaf from a presheaf, known as the \emph{sheafification process} \cite{SR}. This process allows one to use local data encoded in a presheaf to study global properties via the associated sheaf.  
We make frequent use of this notion throughout the article.

 Let $\mathcal{E}$ be a sheaf of abelian groups on a topological space $X$. A section $s$ of $\mathcal{E}$, written as $s \in \mathcal{E}$, means there exists an open set $U \subset X$ such that $s \in \mathcal{E}(U)$. For an open subset $V \subset U$, the restriction of $s$ to $V$ is denoted by $s|_V$, i.e., $s|_V = \operatorname{res}^{\mathcal{E}}_{UV}(s)$, where $\operatorname{res}^{\mathcal{E}}_{UV} : \mathcal{E}(U) \to \mathcal{E}(V)$ is the restriction map. The stalk of $\mathcal{E}$ at a point $x \in X$ is denoted by $\mathcal{E}_x$.  As is customary, we primarily work at the level of spaces of sections.

Let $\mathbb{K}$ be a field of characteristic zero. The constant sheaf on  $X$ with stalks $\mathbb{K}$, is denoted by $\mathbb{K}_X$. A sheaf of $\mathbb{K}$-(Lie) algebras over $X$ is simply denoted by $\mathbb{K}_X$-(Lie) algebra. Similar terminology applies for modules.
   
\subsection{Lie algebroids over ringed spaces} \label{Lie algebroids}
	 Let $\mathcal{O}$ be a sheaf of commutative $\mathbb{K}$-algebras over a topological space $X$. The pair $(X, \mathcal{O})$ is said to be a (commutative) ringed space. The notion of Lie algebroids over (commutative) ringed spaces has been studied in \cite{CV, CRV, MK, AA}.

\label{Smooth spaces}  
We refer to a ringed space \( (X, \mathcal{O}_X) \) as \emph{non-singular} or \emph{smooth} if \( X \) is a smooth manifold, a complex manifold, a smooth algebraic variety, or a smooth scheme of finite type, and \( \mathcal{O}_X \) is its structure sheaf. 

The sheaf of derivations $\mathcal{D}er_{\mathbb{K}_X}(\mathcal{O})$ over $(X, \mathcal{O})$ carries the canonical structure of both an $ \mathcal{O}$-module and a $\mathbb{K}_X$-Lie algebra,  where these structures are compatible through the Leibniz rule. 

With these notations, we recall the definition of Lie algebroids over ringed spaces.

\begin{Def}\label{Category LA}  A Lie algebroid over a ringed space $(X,\mathcal{O})$ is a  $\mathbb{K}_X$-Lie algebra $\mathcal{L}$ that is also an $\mathcal{O}$-module and equipped with an action $\mathcal{L}\otimes_{\mathbb{K}_X} \mathcal{O} \rightarrow \mathcal{O}$ that induces a homomorphism  $$\mathfrak{a}:\mathcal{L}\rightarrow \mathcal{D}er_{\mathbb{K}_X}(\mathcal{O})$$
of $\mathcal{O}$-modules and $\mathbb{K}_X$-Lie algebras. The map $\mathfrak{a}$ is called `the anchor map' and it satisfies the Leibniz rule $$[D,f D^\prime] = f[D,D^\prime]+\mathfrak{a}(D)(f) D^\prime,\quad\text{ for any } f \in \mathcal{O} \text{ and } D,D^\prime \in \mathcal{L}.$$
Fixing the topological space $X$, we denote the Lie algebroid $\mathcal{L}$ $($i.e. the data $(\mathcal{L}, [\cdot, \cdot], \mathfrak{a}))$ over the ringed space $(X,\mathcal{O})$ by the pair $(\mathcal{O}, \mathcal{L})$. Equivalently, $(\mathcal{O}, \mathcal{L})$ can be viewed as a sheaf of Lie-Rinehart algebras over $X$.

A homomorphism of Lie algebroids $$(\phi, \psi) :  (\mathcal{O}_1, \mathcal{L}_1)\rightarrow (\mathcal{O}_2,\mathcal{L}_2)$$  is a sheaf homomorphism of Lie-Rinehart algebras (see \cite{JH} for the local description), i.e. 
	 \begin{itemize}
\item $\psi: \mathcal{L}_1 \rightarrow \mathcal{L}_2$ is a morphism of $\mathcal{O}_1$-modules where $\mathcal{O}_1$ acts on $\mathcal{L}_2$ via $\phi$,
\item $\psi :  (\mathcal{L}_1,[\cdot,\cdot]_1)\rightarrow (\mathcal{L}_2,[\cdot,\cdot]_2)$ is a $\mathbb{K}_X$-Lie algebra homomorphism,
\item compatibility:  $\mathfrak{a}_2 (\psi (D)) (\phi(f))= \phi(\mathfrak{a}_1(D)(f))$, for $f \in \mathcal{O}_1$ and $D \in \mathcal{L}_1$.
\end{itemize}
For a fixed topological space $X$, the Lie algebroids over the ringed spaces associated to $X$ and their morphisms form a category, which we denote by $\mathcal{LA}_X$.
\end{Def} 
\begin{Rem} \label{usual morphisms}
 To follow the usual category of Lie algebroids over a ringed space $(X, \mathcal{O})$ from the existing literature, we consider the following. Form the subcategory 
 \( \mathcal{LA}_{(X, \mathcal{O})} \) of the category $\mathcal{LA}_X$, whose morphisms are taken such that the identity map is the only map on the underlying morphism from $\mathcal{O}$ to $\mathcal{O}$. In particular, the anchor map for an object of \( \mathcal{LA}_{(X, \mathcal{O})} \) forms such a Lie algebroid homomorphism.
\end{Rem}

\emph{Locally free Lie algebroids.}  A Lie algebroid over $(X, \mathcal{O})$ is called a locally free Lie algebroid (of finite rank) if its underlying $\mathcal{O}$-module is locally free (of finite rank). Similarly, we refer to a Lie algebroid as (quasi)coherent Lie algebroid if its underlying $\mathcal{O}$-module structure is (quasi)coherent (see \cite{MK, Jong-Johan-Max-Shin}).

For a (smooth or holomorphic) vector bundle $E$ over a (smooth or complex) manifold $X$, the sheaf of sections of $E$ over $X$ is denoted by $\Gamma_X(E)$. Here, $\Gamma_X$ represents the sheaf of sections functor.
\begin{Exm} The standard Lie algebroid structure on the sheaf $\mathcal{D}er_{\mathbb{K}_X}(\mathcal{O})$ over a ringed space $(X, \mathcal{O})$ is described by the tuple $(\mathcal{D}er_{\mathbb{K}_X}(\mathcal{O}), [\cdot, \cdot]_c, id_X)$, where $[\cdot, \cdot]_c$ is the commutator Lie bracket and the anchor map is  the identity map $id_X$. In particular, if $X$ is a real smooth manifold $($or a complex manifold$)$ and $\mathcal{O}_X$ is its structure sheaf, the sheaf of smooth $($or holomorphic$)$ vector fields $\Gamma_{X}(TX)$ on $X$ is isomorphic to  $\mathcal{D}er_{\mathbb{K}_X}(\mathcal{O}_X)$, where $TX$ is the tangent bundle of $X$ and the base field $\mathbb{K}$ is either $\mathbb{R}$ or $\mathbb{C}$ respectively.
\end{Exm}	
The tangent sheaf 
$\mathcal{T}_X := \left( \mathcal{D}er_{\mathbb{K}_X}(\mathcal{O}_X), [\cdot,\cdot]_c, \mathrm{id}_X \right)$
forms a locally free Lie algebroid of finite rank over a non-singular or smooth space $(X, \mathcal{O}_X)$. It is an algebraic generalization of the notion $\Gamma_{X}(TX)$.

Now, we recall important examples to show that the Definition \ref{Category LA} unifies the definitions of Lie algebroids over smooth manifolds, complex manifolds, analytic spaces, algebraic varieties, and schemes. 
 \begin{Exm}[Smooth Lie algebroids] \label{Classical cases}
           Lie algebroids over a smooth manifold $X$ is equivalent to locally free Lie algebroids of finite rank over the ringed space $(X, C_X^\infty)$, where $ C^\infty_X$ is the sheaf of  $C^\infty$-functions on $X$.
\end{Exm}

\begin{Exm}[Holomorphic Lie algebroids]
  Let $(X,\mathcal{O}_X)$ be a complex manifold equipped with the structure sheaf $\mathcal{O}_X$ that is the sheaf of holomorphic functions on $X$. A locally free Lie algebroid of finite rank over $(X,\mathcal{O}_X)$ is precisely a holomorphic Lie algebroid, we refer to \cite{SC, LG-PX, PT, BP} for more details.  
\end{Exm}
  Consider a holomorphic Lie algebroid $(\mathcal{O}_X,\mathcal{L})$  over a complex manifold $X$. Using the underlying smooth structure, one have a smooth Lie algebroid $(C^{\infty}_{X, \mathbb{C}}, \mathcal{L}^{\infty})$, where $C^{\infty}_{X, \mathbb{C}}$  is the sheaf of $\mathbb{C}$-valued smooth functions on $X$ and 
	$\mathcal{L}^{\infty}:= C^{\infty}_{X, \mathbb{C}} \otimes_{\mathcal{O}_X} \mathcal{L}$	is the associated 	$C^{\infty}_X$-module (see \cite{LG-PX}). Thus, we have canonical monomorphism
    $$(\mathcal{O}_X,\mathcal{L}) \hookrightarrow (C^{\infty}_{X, \mathbb{C}}, \mathcal{L}^{\infty}).$$

 \begin{Exm}[Poisson manifolds] \label{Cotantgent Poisson}
		The cotangent sheaf $\Omega^1_X$ $($i.e., the sheaf of differential $1$-forms$)$ over a (smooth or holomorphic) Poisson manifold $X$ with Poisson bi-vector field $\pi \in H^0(X, \wedge^2 \mathcal{T}_X)$, has a canonical Lie algebroid structure.
        The anchor map $\tilde{\pi}: \Omega^1_X \rightarrow \mathcal{T}_X$ and the Lie bracket on $\Omega^1_X$ is defined by $$\tilde{\pi}(\omega) (f)= \pi(\omega, df),$$   $$[\omega, \omega']:= L_{\tilde{\pi}(\omega)} \omega' - L_{\tilde{\pi}(\omega')} \omega - d(\pi(\omega, \omega')),$$ where  $ f \in \mathcal{O}_X$ and $\omega, \omega^\prime \in \Omega^1_X.$
        We refer to  \cite{RF, LG-PX, BP} for more details.
	\end{Exm}	

Above examples are locally free Lie algebroids over the underlying ringed spaces. Now, we discuss an example, which is not necessarily locally free over some spaces having possible singularities.  			
            
    \begin{Exm}[Analytic spaces and Singular Foliations] \label{Foliation}
    Let $(X,\mathcal{O}_X)$ be a complex manifold. The zero locus $Y:= V(\mathcal{I})$ of an coherent sheaf of ideals $\mathcal{I} \subset \mathcal{O}_X$ is not necessarily a submanifold, and may have singularities as well. The pair $(Y, \mathcal{O}_Y)$ forms an analytic subspace of $(X,\mathcal{O}_X)$,  where the quotient sheaf $\mathcal{O}_Y := \mathcal{O}_X/{\mathcal{I}}$ considered to be the structure sheaf on $Y$ \cite{CW, BP}. The tangent sheaf $\mathcal{T}_Y := \mathcal{D}er_{\mathbb{C}_Y}(\mathcal{O}_Y)$ forms a Lie algebroid over $(Y, \mathcal{O}_Y)$. Consider the sheaf of logarithmic derivations \cite{Francisco, BP, AA} given by 
				$$\mathcal{T}_X(-log~ Y) := \{D \in \mathcal{T}_X : D(\mathcal{I}) \subset \mathcal{I}\} \hookrightarrow \mathcal{T}_X$$
    with the canonical Lie algebroid structure. It is associated with $\mathcal{T}_Y$ via the Lie algebroid epimorphism 
	$$(\rho_0, \rho_1) : (\mathcal{O}_X, \mathcal{T}_X(-log~ Y)) \rightarrow (\mathcal{O}_Y, \mathcal{T}_Y)$$ 
	defined by $\rho_1 (D) = \tilde{D}$, for $D \in\mathcal{T}_X(-log~ Y)$ as $\tilde{D}(\bar{f}) = \overline{D(f)}$, where $\rho_0(f)=\bar{f}$ for $f \in \mathcal{O}_X$, is the quotient map.
    Notice that, both the sheaf \((\mathcal{O}_X, \mathcal{T}_X(-\log Y) )\) and \((\mathcal{O}_Y, \mathcal{T}_Y) \) form coherent Lie algebroids.

    The sheaf of logarithmic derivations is an example of a generalized involutive distribution or singular foliation (see \cite{LG, BP, JV}). Recall that a singular foliation  $\mathcal{F}$ on a (smooth or complex) manifold $(X, \mathcal{O}_X)$ is an $\mathcal{O}_X$-submodule of the Lie algebroid $\mathcal{T}_X$  $($or, $\Gamma_{X}(TX))$ that is stable under the Lie bracket and locally finitely generated. Thus, any singular foliation on a (smooth or complex) manifold $(X, \mathcal{O}_X)$ forms a (coherent) Lie subalgebroid of $\mathcal{T}_X$. We refer to $\mathcal{T}_X(-\log~ Y)$ as the logarithmic foliation of $X$ along $Y$. 
	\end{Exm}

\begin{Exm}[Lie-Rinehart algebras]\label{LR on affine schemes} Let $(R, L)$ be a Lie-Rinehart algebra for a commutative $\mathbb{K}$-algebra $R$. The associated affine scheme is given by $(X, \mathcal{O}_X)$, where $X=Spec(R)$ with the structure sheaf of regular functions $\mathcal{O}_X$ induced by localization of the ring $R$. If we localize the $R$-module $L$ to obtain the $\mathcal{O}_X$-module $\mathcal{L}$, then a Lie algebroid structure is induced on $\mathcal{L}$ over the ringed space $(X, \mathcal{O}_X)$, see \cite{CW} for more details. Conversely, starting with a Lie algebroid $\mathcal{L}$ over the ringed space $(X,\mathcal{O}_X)$, we get a Lie-Rinehart algebra $\mathcal{L}(U)$ over $\mathcal{O}_X(U)$ for each open set $U\subset X$. With these observations, one can show that the category of Lie-Rinehart algebras over $R$ and the category of quasicoherent Lie algebroids over the associated affine scheme $(X,\mathcal{O}_X)$ are equivalent. 
More generally, using this correspondence locally over affine open covers, Lie algebroids have been studied on Noetherian separated schemes and on schemes of finite type (see \cite{MK, UB}).
\end{Exm}

\subsection{Universal enveloping algebroid of a Lie algebroid} \label{Universal enveloping algebroid}
Similar to the universal enveloping algebra of a Lie algebra, the universal enveloping algebra of a Lie-Rinehart algebra is fundamental in studying its homological algebra (see \cite{GR, JH, HM}). The algebro-geometric counterpart of this notion is known as the universal enveloping algebroid (see \cite{MK, DRV, CV, BP, UB, BRT, AA}).

			Let $(\mathcal{L},[\cdot,\cdot],\mathfrak{a})$ be  a Lie algebroid  over  $(X, \mathcal{O})$.  For each open set $U$ of $X$, we find  the universal enveloping algebra $\mathcal{U}(\mathcal{O}(U),\mathcal{L}(U))$ of the Lie-Rinehart algebra $(\mathcal{O}(U),\mathcal{L}(U))$ (see \cite{GR,JH,MM}). Apart from the well known constructions (see Remark \ref{0-twist uni alg}), the universal enveloping algebra $\mathcal{U}(\mathcal{O}(U),\mathcal{L}(U))$ of the Lie-Rinehart algebra $(\mathcal{O}(U),\mathcal{L}(U))$ can be described as follows (see \cite{NK-PS, TS}). It is given by the quotient algebra 
            $$T_{\mathcal{O}(U)} (\mathcal{L}(U) \otimes_{\mathbb{K}} \mathcal{O}(U))/{\mathcal{J}(U)},$$ where $T_{\mathcal{O}(U)} (\mathcal{L}(U) \otimes_{\mathbb{K}} \mathcal{O}(U))$ is the tensor $\mathcal{O}(U)$-algebra of the $\mathcal{O}(U)$-bimodule $\mathcal{L}(U) \otimes_{\mathbb{K}} \mathcal{O}(U)$, and $\mathcal{J}(U)$ is the two sided ideal generated by the elements
            $$\eta_U(D) \otimes \eta_U(D')- \eta_U(D') \otimes \eta_U(D)- \eta_U([D, D']),~ \eta_U(D) f- f \eta_U(D)- \mathfrak{a}_U(D)(f),$$
            where $f \in \mathcal{O}(U)$, $D, D' \in \mathcal{L}(U)$, ~ the map $\mathfrak{a}_U:\mathcal{L}(U) \rightarrow Der_{\mathbb{K}}(\mathcal{O}(U))$ is the anchor for $(\mathcal{O}(U),\mathcal{L}(U))$ and the map $\eta_U: \mathcal{L}(U) \rightarrow \mathcal{L}(U) \otimes_{\mathbb{K}} \mathcal{O}(U)$ is defined by $D \mapsto D \otimes 1$.

            For a Lie algebroid $(\mathcal{O}, \mathcal{L})$, the sheafification of the cannonical presheaf: 
			$$U \mapsto \mathcal{U}(\mathcal{O}(U),\mathcal{L}(U)),$$ is known as the universal enveloping algebroid of the Lie algebroid $\mathcal{L}$, and denoted by $\mathscr{U}(\mathcal{O}, \mathcal{L})$. The underlying presheaf $\mathcal{U}(\mathcal{O},\mathcal{L})$ satisfies the compatibility condition that for every open set $V\subset U$ there is a canonical restriction map 
$res_{UV}^\mathcal{U}:\mathcal{U}(\mathcal{O}(U),\mathcal{L}(U))\rightarrow \mathcal{U}(\mathcal{O}(V),\mathcal{L}(V))$ defined by 
$$res_{UV}^\mathcal{U}(f~ \overline{D_1} \cdots \overline{D_n} )= f|_V ~\overline{D_1|_V} \cdots \overline{D_n|_V},$$
where $f\in \mathcal{O}_X(U),~ D_i\in \mathcal{L}(U),~ 
\overline{D_i} = \iota_{\mathcal{L}(U)}(D_i),$ for $i=1, \dots, n$
and $\iota_{\mathcal{L}(U)}:\mathcal{L}(U) \rightarrow \mathcal{U}(\mathcal{O}_X(U),\mathcal{L}(U))$ is the canonical map given by $D \mapsto \bar{D}$.

			From the construction of $\mathscr{U}(\mathcal{O}, \mathcal{L})$, it is an associative $\mathbb{K}_X$-algebra and $\mathcal{O}$-bimodule. Moreover,
			there is a canonical  $\mathbb{K}_X$-algebra monomorphism $\iota: \mathcal{O} \hookrightarrow \mathscr{U}(\mathcal{O},\mathcal{L})$ and an $\mathcal{O}$-linear map $\iota_{\mathcal{L}}: \mathcal{L} \rightarrow \mathscr{U}(\mathcal{O},\mathcal{L})$. 
			The associative $\mathbb{K}_X$-algebra $\mathscr{U}(\mathcal{O}_X, \mathcal{L})$ is generated by $\mathcal{O}$ and $\iota_{\mathcal{L}}(\mathcal{L})$ satisfy the following identities:
            \begin{align} \label{relations in U}
             \bar{D} ~ \bar{D'} - \bar{D'} ~ \bar{D}= \overline{[D,~D']},~  \bar{D} ~ f - f ~ \bar{D}= \mathfrak{a}(D)(f),   
            \end{align}
			where $D, D' \in \mathcal{L}$, ~$f \in \mathcal{O}$, and $\bar{D}= \iota_{\mathcal{L}}(D)$ for all $D\in \mathcal{L}$.
			
			Hence, the map $\iota_{\mathcal{L}}$ can also be viewed as a $\mathbb{K}_X$-Lie algebra homomorphism.
			\begin{Rem}\label{Uni prop}
				The universal enveloping algebroid $\mathscr{U}(\mathcal{O}, \mathcal{L})$ of a Lie algebroid $(\mathcal{O}, \mathcal{L})$ is characterized by the following universal property: \label{Universal}
				Let $\mathcal{A}$ be a unital associative $\mathbb{K}_X$-algebra with sheaf homomorphisms $\phi: \mathcal{O} \rightarrow \mathcal{A}$ of $\mathbb{K}_X$-unital algebras and $\psi: (\mathcal{L}, [\cdot,\cdot]) \rightarrow (\mathcal{A}, [\cdot,\cdot]_c)$ of  $\mathbb{K}_X$-Lie algebras such that $\phi(f)\psi(D) = \psi(fD)$ and $[\psi(D), \phi(f)]_c = \phi(\mathfrak{a}(D)(f))$ holds for $f \in \mathcal{O}$ and $D \in \mathcal{L}$. Then, there exists a unique  homomorphism of unital $\mathbb{K}_X$-algebras $\widetilde{\psi} : \mathscr{U}(\mathcal{O}, \mathcal{L}) \rightarrow \mathcal{A}$ such that $\widetilde{\psi} \circ \iota = \phi$ and $\widetilde{\psi} \circ \iota_{\mathcal{L}} = \psi$.
			\end{Rem}
			 
			The $\mathbb{K}_X$-algebra $\mathscr{U}(\mathcal{O},\mathcal{L})$  has a natural filtration of $\mathcal{O}$-modules: 
            \begin{align} \label{natural filtration}
			    \mathcal{O} = \mathscr{U}_{0}(\mathcal{O},\mathcal{L}) \subset \mathscr{U}_{1}(\mathcal{O},\mathcal{L}) \subset \mathscr{U}_{2}(\mathcal{O},\mathcal{L}) \subset \cdots,
			\end{align}
			where $\mathscr{U}_{n}(\mathcal{O},\mathcal{L})$ is spanned by the sections of $\mathcal{O}$ and the powers $\iota_{\mathcal{L}}(\mathcal{L})^m$ for $m = 1,2, \dots ,n$.
			The  direct sum of the quotient sheaves associated with each of the consecutive sheaves appear in the above filtration of $\mathscr{U}(\mathcal{O}_X,\mathcal{L})$ forms a sheaf of graded algebras. It is denoted by $gr(\mathscr{U}(\mathcal{O}_X,\mathcal{L}))$. Therefore,
				$$gr(\mathscr{U}(\mathcal{O},\mathcal{L})) = \bigoplus_{n\geq 0} \mathscr{U}_{n}(\mathcal{O},\mathcal{L})/\mathscr{U}_{n-1}(\mathcal{O},\mathcal{L}),$$
				where $\mathscr{U}_{-1}(\mathcal{O},\mathcal{L})$ is the zero sheaf by convention. It is a commutative associative untial $\mathcal{O}$-algebra.
			
			\begin{Exm}
				For a non-singular ringed space $(X,\mathcal{O}_X)$, the universal enveloping algebroid of $\mathcal{T}_X$ is isomorphic to the sheaf of differential operators $\mathcal{D}_X$ on $X$ (i.e. the sheaf of differential operators over $\mathcal{O}_X$, sometimes denoted as $\mathcal{D}iff(\mathcal{O}_X))$ \cite{MK, TS, SR}, equipped with a canonical isomorphism $$\mathscr{U}(\mathcal{O}_X, \mathcal{T}_X) \cong \mathcal{D}iff(\mathcal{O}_X)=:\mathcal{D}_X.$$

   In particular, let $X$ be an $n$-dimensional $C^\infty$-manifold or a Stein manifold. Then $L:= Der_{\mathbb{K}}(\mathcal{O}_X(X))$ is the $(\mathbb{K}, \mathcal{O}_X(X))$-Lie-Rinehart algebra of first order homogeneous differential operators on $X$ and its universal enveloping algebra $\mathcal{U}(\mathcal{O}_X(X), L)$ is isomorphic to the $\mathbb{K}$-algebra (and $\mathcal{O}_X(X)$-module) of differential operators on $X$. In a local coordinate system $(U, (x_1, \dots ,x_n))$, elements of $\mathcal{U}_{k}(\mathcal{O}_X(U), Der_{\mathbb{K}}(\mathcal{O}_X(U)))$ are of the form
	$$\displaystyle \sum_{i_1 + \dots +i_n \leq k}f_{i_1 \dots i_n}(\partial/\partial x_1)^{i_1} \circ \cdots \circ (\partial/\partial x_n)^{i_n}, ~~\text{~~for some}~~ f_{i_1 \dots i_n} \in \mathcal{O}_X(U).$$
		Moreover, this description holds locally for $\mathcal{D}_X$ over a complex manifold $X$ via Stein open covers \cite{BRT}.
			\end{Exm}
   \begin{Exm} We give a description for $\mathscr{U}(\mathcal{O}_X, \mathcal{T}_X)$ or $\mathcal{D}_X$ over a smooth affine scheme $(X, \mathcal{O}_X)$, using the algebraic approach described in Example \ref{LR on affine schemes}, in the following.
      An algebraic generalizations of differential operators, the algebra of differential operators \( \mathrm{Diff}_{\mathbb{K}}(R) \) for a commutative \( \mathbb{K} \)-algebra \( R \), is also denoted by \( \mathcal{D}_{R/\mathbb{K}} \), introduced by Grothendieck (see \cite{LM},\cite{TS}). In this algebra, the following relations hold for any  \( f \in R \):
$$[D, f]_c = D(f) \in R ~~\text{for}~~D \in {Der}_{\mathbb{K}}(R), 
~~~~~\text{and}~~ [\tilde{D}, f]_c \in  \mathcal{D}^{\leq n}_{R/\mathbb{K}}~~~~~\text{~~~~for all}~~\tilde{D} \in \mathcal{D}^{\leq n+1}_{R/\mathbb{K}},$$
where \( [\cdot,\cdot]_c \) denotes the commutator and  $(\mathcal{D}_{R/\mathbb{K}}, \mathcal{D}^{\leq n}_{R/\mathbb{K}})$ is the canonical filtered algebra structure.
Recall that under certain algebraic conditions on \( R \), the canonical \( \mathbb{K} \)-algebra epimorphism
\[
\mathrm{Diff}_{\mathbb{K}}(R) \rightarrow \mathcal{U}(R, \mathrm{Der}_{\mathbb{K}}(R))
\]
becomes an isomorphism. This holds, for instance, when \( R = C^\infty(M) \) for a smooth manifold \( M \), or when \( R = \mathbb{K}[x_1, \dots, x_n] \), or more generally, when \( R \) is a regular (or a smooth) \( \mathbb{K} \)-algebra (see \cite{HM, TS}).

   \end{Exm}
			\begin{Exm} \label{sheaf of log diff operators}
				The sheaf of logarithmic derivations and sheaf of logarithmic differential operators (represents differential operators along the foliation) are denoted as $\mathcal{T}_X(-log~ Y)$ and $\mathcal{D}_X(-log~ Y)$ respectively for a principal divisor $Y$ in some complex manifold or smooth algebraic variety $X$ \cite{Francisco, AA}.
				
				In the case of a free divisor $Y$ in $X$ (i.e. $\mathcal{T}_X(-log~ Y)$ is locally free $\mathcal{O}_X$-module \cite{Francisco, CMD,BP,AA})  we have (sheafifying the local description given for the module of logarithmic derivations in \cite{LM}) the isomorphism
				$$\mathscr{U}(\mathcal{O}_X, \mathcal{T}_X(-log~ Y)) \cong  \mathcal{D}_X(-log~ Y).$$
				
			\end{Exm}

			\begin{Exm}\label{free Lie algebroid} In \cite{MK}, the notion of path algebroid $\mathcal{P}_X$ of a smooth manifold or smooth algebraic variety $X$ is constructed as the free Lie algebroid generated by tangent sheaf $\mathcal{T}_X$.  It forms a locally free (quasicoherent) Lie algebroid over $(X,\mathcal{O}_X)$, possibly of infinite rank. The universal enveloping algebroid $\mathscr{U}(\mathcal{O}_X,\mathcal{P}_X) =:\mathbb{D}_X$ is described as sheaf of non-commutative differential operators on $X$.
				
			\end{Exm} 
            Recall the Poincaré-Birkhoff-Witt (Generalized PBW)  theorem for Lie algebroids over a ringed space \cite{AA}.
    If $(\mathcal{O}, \mathcal{L})$ is a locally free Lie algebroid, then the graded algebra $gr(\mathscr{U}(\mathcal{O}, \mathcal{L}))$ is isomorphic to the symmetric algebra 	$\mathscr{S}_{\mathcal{O}}\mathcal{L}$ as graded $\mathcal{O}$-algebras, given by the symmetrization map $($known as the PBW map$)$
	\begin{align}  \label{PBW}
		 \psi_{\mathcal{L}}:\mathscr{S}_{\mathcal{O}}\mathcal{L} \rightarrow gr(\mathscr{U}(\mathcal{O}, \mathcal{L})
					\end{align}
					\begin{center}
					    $D_1 \otimes \cdots \otimes D_k \mapsto  \frac{1}{k!} \displaystyle\sum_{\sigma \in S_k} \overline{D_{\sigma(1)}} \cdots \overline{D_{\sigma(k)}}$
					\end{center}	
					where $D_1, \dots , D_k$ are sections of $\mathcal{L}$ and $\overline{D_{i}}$ is the associated class of $D_i$ for $i= 1, \dots, k$.
				\begin{Rem}\label{embed}
					For a locally free Lie algebroid $\mathcal{L}$ over $(X, \mathcal{O})$, the homomorphism $\iota_{\mathcal{L}}: \mathcal{L} \rightarrow \mathscr{U}(\mathcal{O},\mathcal{L})$ become an embedding $($as $\mathbb{K}_X$-Lie algebras and $\mathcal{O}$-modules$)$. 
				\end{Rem}
				\begin{Rem}
						If $\mathcal{L}$ is locally free $\mathcal{O}$-module then by comparing the $\mathcal{O}$-modules in degree $1$ of the PBW isomorphism we get  $\mathcal{L} \cong \mathscr{U}_{1}(\mathcal{O}, \mathcal{L})/{\mathcal{O}}$. Thus, we have the short exact sequence  $$0 \rightarrow \mathcal{O} \rightarrow \mathscr{U}_{1}(\mathcal{O}, \mathcal{L}) \rightarrow \mathcal{L} \rightarrow 0$$ of $\mathcal{O}$-modules, which is also a
					split exact sequence, i.e.  $\mathscr{U}_{1}(\mathcal{O}, \mathcal{L}) \cong \mathcal{O} \oplus \mathcal{L}$. From here we get a canonical embedding of $\mathcal{L}$ in $\mathscr{U}(\mathcal{O}, \mathcal{L})$ and the $\mathcal{O}$-module $\mathcal{O} \oplus \mathcal{L}$ has a canonical Lie algebroid structure given as 
                    $$[(f, D), (g, D')]_c:= (\mathfrak{a}({D})(g)- \mathfrak{a}({D}')(f), [D,D']) ~~\text{~~for}~~f, g \in \mathcal{O},~ D, D' \in \mathcal{L}.$$
				\end{Rem}
 \subsection{Chevalley-Eilenberg-de Rham hypercohomology}
To consider Lie algebroid cohomology for a Lie algebroid $(\mathcal{O}, \mathcal{L})$ with coefficient in some $\mathcal{O}$-module $\mathcal{E}$, first we consider the notion of the \emph{Atiyah algebroid} of $\mathcal{E}$. Then we consider an algebro-geometric analogue of the standard cochain complex of smooth Lie algebroids (see \cite{KM}) and consider its hypercohomology (see \cite{VID}). This is refereed as the \emph{Chevalley-Eilenberg-de Rham hypercohomology} (see \cite{UB, BRT, BP, PT}).

\subsubsection{Representations of a Lie algebroid}
For a (quasicoherent) $\mathcal{O}$-module $\mathcal{E}$, we form a Lie algebroid consisting of the sheaf of differential operators on $\mathcal{E}$ of order $\leq 1$ having scalar symbols \cite{MK,UB}, i.e. 
$$\mathcal{A}t(\mathcal{E}):= \{D\in \mathcal{E}nd_{\mathbb{K}_X}(\mathcal{E})~|~D(fs)=fD(s)+\sigma_D(f)s~~\mbox{ for a unique}~~  \sigma_D\in \mathcal{D}er_{\mathbb{K}_X}(\mathcal{O}),~~\mbox{where}~~ f \in \mathcal{O},~s\in \mathcal{E}\},$$
(thus, here $\sigma_D(f)=[D, f]_c \in \mathcal{O}$ for $D \in \mathcal{A}t(\mathcal{E})$,~ $f \in \mathcal{O}$ holds) 
with the anchor map defined by 
$$\sigma: \mathcal{A}t(\mathcal{E}) \rightarrow \mathcal{D}er_{\mathbb{K}_X}(\mathcal{O})~~\text{where}~~D \mapsto \sigma_D$$
and the Lie bracket is commutator bracket. This Lie algebroid structure is so-called Atiyah algebroid of the $\mathcal{O}$-module $\mathcal{E}$. In smooth vector bundles or principal bundles context, one can find such notions in \cite{KM-KS}.

For a Lie algebroid $(\mathcal{L}, [\cdot, \cdot], \mathfrak{a})$ over $(X, \mathcal{O})$, an $\mathcal{L}$-connection on $\mathcal{E}$ is defined by an $\mathcal{O}$-linear map (see \cite{MK,UB})
\begin{align}\label{L-connection}
	\nabla: \mathcal{L}\rightarrow \mathcal{A}t(\mathcal{E}) ~~\text{~~sends}~~D\mapsto \nabla_D,
\end{align}
satisfying the Leibniz rule $$\nabla_D(f~s)= f~\nabla_D(s)+\mathfrak{a}(D)(f)~s,$$ for  sections $f \in \mathcal{O}$, $D \in \mathcal{L}$ and $s \in \mathcal{E}$ $($to ensure compatibility, we choose $\sigma \circ \nabla= \mathfrak{a}$, i.e. $\sigma_{\nabla_{D}}=\mathfrak{a}(D)$ for $D\in \mathcal{L})$. Equivalently, it is described by a $\mathbb{K}_X$-linear map, if in addition $\mathcal{L}$ is a locally free $\mathcal{O}$-module of finite rank, as follows (see \cite{BP})
\begin{align} \label{L-connection special}
	d_{\nabla} : \mathcal{E} \rightarrow \Omega_{\mathcal{L}}^1 \otimes_{\mathcal{O}}\mathcal{E}
\end{align}
satisfying the  Leibniz rule 
$d_{\nabla}(f~s)=f~d_{\nabla} (s) + {\mathfrak{a}^*}(df)\otimes s,$ where $\Omega^1_{\mathcal{L}}:=\mathscr{H}om_{\mathcal{O}}(\mathcal{L}, \mathcal{O})$ $($or, $\mathcal{L}^*)$ is the dual of a Lie algebroid $\mathcal{L}$, together with  $\mathcal{O}$-module homomorphism ${\mathfrak{a}^*}$, the dual of the anchor map $\mathfrak{a}$. The association is given by $d_{\nabla}(s)(D)=\nabla_D(s)$, for  $f \in \mathcal{O}$, $D \in \mathcal{L}$.
The curvature of the $\mathcal{L}$-connection $\nabla$ is the $\mathcal{O}$-linear map $$R_{\nabla}:\wedge^2_{\mathcal{O}} \mathcal{L} \rightarrow \mathcal{E}nd_{\mathbb{K}_X}(\mathcal{E})$$ defined as 
\begin{align} \label{curvature}
  R_{\nabla}(D\wedge D')= [\nabla_D,\nabla_{D'}]_c - \nabla_{[D,D']_c},  
\end{align}
 for any two sections $D,D'$ of $\mathcal{L}$. One verifies that the curvature $R_{\nabla} \in \mathscr{H}om_{\mathcal{O}}(\wedge^2_{\mathcal{O}}\mathcal{L}, \mathcal{E}nd_{\mathcal{O}}(\mathcal{E}))$ (see \cite{HM-Chernclass}).  

An $\mathcal{L}$-connection $\nabla$ on $\mathcal{E}$ is said to be flat if the map (\ref{L-connection}) is a Lie algebroid homomorphism  (i.e. the $\mathcal{L}$-curvature is zero), i.e. the map (\ref{L-connection}) additionally satisfies
\label{L-module} $R_{\nabla} (D, D')=0$, for all $D, D' \in \mathcal{L}$.
In this case, $(\mathcal{E},\nabla)$ is referred to as a representation of $\mathcal{L}$ or, equivalently, an $\mathcal{L}$-module.

If $(\mathcal{E}, \nabla)$ is an $\mathcal{L}$-module, then the Lie algebroid morphism
$\nabla : \mathcal{L} \rightarrow \mathcal{A}t(\mathcal{E}) $ 
extends to an $\mathcal{O}$-linear  homomorphism of $\mathbb{K}_X$-algebras (using Remark \ref{Universal})
$$\tilde{\nabla} : \mathscr{U}(\mathcal{O},\mathcal{L}) \rightarrow \mathscr{E}nd_{\mathbb{K}_X}(\mathcal{E}),$$
making $\mathcal{E}$ into a left $\mathscr{U}(\mathcal{O},\mathcal{L})$-module. If $\mathcal{L}$ is locally free $\mathcal{O}$-module and $\mathcal{E}$ be a $\mathscr{U}(\mathcal{O},\mathcal{L})$-module, the restriction of the action of $\mathscr{U}(\mathcal{O},\mathcal{L})$ to $\mathcal{L}$ (using the canonical embedding of $\mathcal{L}$ in $\mathscr{U}(\mathcal{O},\mathcal{L})$ as described in  Remark \ref{embed}), 
provides an $\mathcal{L}$-module structure on $\mathcal{E}$. In this case, the category of $\mathcal{L}$-modules and the category of left $\mathscr{U}(\mathcal{O},\mathcal{L})$-modules are equivalent. It helps in studying  homological algebra with $\mathcal{L}$-modules ($\mathcal{L}$ is a sheaf of non-associative $\mathbb{K}$-algebras) by treating them as modules over the associative $\mathbb{K}_X$-algebra $\mathscr{U}(\mathcal{O},\mathcal{L})$.

\subsubsection{Chevalley-Eilenberg-de Rham complex}
For a Lie algebroid $\mathcal{L}:=(\mathcal{L}, [\cdot, \cdot], \mathfrak{a})$ over $(X,\mathcal{O})$ with a representation $(\mathcal{E},\nabla)$, consider the cochain complex (consists with  $\mathcal{O}$-modules with a $\mathbb{K}_X$-linear degree $1$ map), a generalization of the well known \emph{Chevalley-Eilenberg-de Rham complex} \cite{KM} is given as follows.

The \emph{Chevalley-Eilenberg-de Rham complex} of $(\mathcal{O}, \mathcal{L})$ with coefficient in $\mathcal{E}$ is 
\begin{align}\label{Chevalley-Eilenberg-de Rham complex}
   \Omega_{\mathcal{L}}^\bullet(\mathcal{E}):=(\wedge^\bullet_{\mathcal{O}}\mathcal{L}^* \otimes_{\mathcal{O}} \mathcal{E}, d_{\mathcal{L}}), 
\end{align}
where the differential $d_{\mathcal{L}}: \wedge^\bullet_{\mathcal{O}}\mathcal{L}^* \otimes_{\mathcal{O}} \mathcal{E} \rightarrow \wedge^{\bullet +1}_{\mathcal{O}}\mathcal{L}^* \otimes_{\mathcal{O}} \mathcal{E}$  is given by
\begin{align*}
	\begin{split}
		d_{\mathcal{L}}(\omega)(D_1\wedge \cdots \wedge D_{k+1}) & = \sum^{k+1}_{i=1}(-1)^{i+1}~\nabla_{D_i}(\omega(D_1\wedge \cdots \wedge \hat{D_i} \wedge \cdots \wedge D_{k+1}))\\
		&+ \sum_{i< j}(-1)^{i+j}~\omega([D_i, D_j]\wedge D_1\wedge \cdots \wedge \hat{D_i}\wedge \cdots \wedge \hat{D_j}\wedge \cdots \wedge D_{k+1}),
	\end{split}
\end{align*}
where $D_1, \dots, D_{k+1} \in \mathcal{L}$ and
$\omega \in \wedge^{k}_{\mathcal{O}} \mathcal{L}^*$, 
and $\nabla$ is the flat $\mathcal{L}$-connection on $\mathcal{E}$ (see the notion in (\ref{L-connection})).  

Notice that, the differential $d_{\mathcal{L}}$ is a $\mathbb{K}_X$-linear map but not an $\mathcal{O}$-module homomorphism, satisfies the graded Leibniz rule: for $\omega_1 \in \wedge^{k}_{\mathcal{O}} \mathcal{L}^*$ and $\omega_2 \in  \wedge^{l}_{\mathcal{O}} \mathcal{L}^*$, where $k, l \in \mathbb{N} \cup \{0\}$,
$$d_{\mathcal{L}}(\omega_1 \wedge \omega_2)= d \omega_1 \wedge \omega_2 + (-1)^k \omega_1 \wedge d\omega_2.$$
Thus, the cochain complex with the exterior product $(\Omega_{\mathcal{L}}^\bullet(\mathcal{E}), \wedge)$ forms a differential graded algebra.

The associated hypercohomology  of the cochain complex (\ref{Chevalley-Eilenberg-de Rham complex}) of $\mathcal{O}$-modules is called the Chevalley-Eilenberg-de Rham cohomology or Lie algebroid cohomology of $\mathcal{L}$ with coefficient in $\mathcal{E}$, and denoted by $\mathbb{H}^\bullet(\mathcal{L}, \mathcal{E})$. In particular, we get the Chevalley-Eilenberg cohomology when $\mathcal{L}$ is a Lie algebroid over a point; and the de Rham cohomology  when  $\mathcal{L}=\mathcal{T}_X$ and $\mathcal{E}=\mathcal{O}_X$ over a non-singular space $X$ $($see Section \ref{Smooth spaces}$)$.

Let $(\mathcal{E}, \nabla)$ be an $\mathcal{L}$-connection and $\omega \in \mathcal{Z}^2(\mathcal{L}, \mathcal{O}):=\mathcal{K}er(d_{\mathcal{L}}: \wedge^2_{\mathcal{O}}\mathcal{L}^*  \rightarrow \wedge^{3}_{\mathcal{O}}\mathcal{L}^*)$.  If the following holds
\begin{align} \label{curvature type}
  R_{\nabla}(D \wedge D')(s)= \omega(D, D')~s,
\end{align}
 $($see $(\ref{curvature}))$ for all $D, D' \in \mathcal{L}$ and $s \in \mathcal{E}$, we say the $\mathcal{L}$-connection $\nabla$ on $\mathcal{E}$ is of \emph{curvature type} $\omega$.

\subsubsection{Abelian Lie algebroid extensions.}
We denote the set of equivalence classes of abelian extensions of a Lie algebroid $\mathcal{L}$ by a flat $\mathcal{L}$-connection $(\mathcal{E}, \nabla)$ as $Ext^1(\mathcal{L}, \mathcal{E},\nabla)$. A representative  $\mathcal{L}'$ of $ Ext^1(\mathcal{L}, \mathcal{E},\nabla)$ corresponds to an abelian extension given by a short exact sequence (s.e.s.) of Lie algebroids
$$0 \rightarrow \mathcal{E} \rightarrow \mathcal{L}'\rightarrow \mathcal{L} \rightarrow 0$$
up to isomorphism in the usual sense, where $\mathcal{E}$ is the trivial Lie algebroid. Thus, as a local description, on the level of stalks, for each point $x \in X$, we have the s.e.s. of $(\mathbb{K}, \mathcal{O}_{X,x})$-Lie-Rinehart algebras
$$0 \rightarrow \mathcal{E}_x \rightarrow {\mathcal{L}'}_x\rightarrow \mathcal{L}_x \rightarrow 0.$$
We recall the following result about abelian extensions of $\mathcal {L}$ by $\mathcal{E}$ (see \cite{AS, PT, HM}).
Let $(\mathcal{O}, \mathcal{L})$ be a locally free Lie algebroid   with an $\mathcal{L}$-module $(\mathcal{E}, \nabla)$. Then, we get a canonical isomorphism of $\mathbb{K}$-vector spaces
\begin{align} \label{Ext-H^2}
   \mathbb{H}^2(\mathcal{L}, \mathcal{E})\cong Ext^1(\mathcal{L}, \mathcal{E},\nabla). 
\end{align}

			\section{Sheaf of quantum Poisson algebras} \label{Sec 3}
			In this section, we consider sheaves of quantum Poisson algebras \cite{GP}, which we refer to as `quantum Poisson algebroids'. As a first example, we examine the canonical quantum Poisson algebroid structure on the universal enveloping algebroid of a Lie algebroid. As another crucial example, we describe the canonical quantum Poisson algebroid structure on the sheaf of algebras of differential operators on a sheaf of modules. Finally, we demonstrate a categorical correspondence between the category of Lie algebroids and the category of quantum Poisson algebroids. We briefly describe the case for smooth, holomorphic, and algebraic contexts. 

\begin{Def} \label{Category QP} A quantum Poisson algebroid $(\mathcal{D},\mathcal{D}^i)$ is a sheaf of quantum Poisson algebras over a topological space $X$. That is,
				$\mathcal{D}=\cup ^{\infty}_{i=0}\mathcal{D}^i$ is an associative $\mathbb{K}_X$-algebra with the filtration
				$$\mathcal{D}^0 \subset \mathcal{D}^1 \subset \cdots \subset \mathcal{D}^n \subset \cdots \subset \mathcal{D}$$
				such that the following properties hold
                $$\mathcal{D}^i \cdot \mathcal{D}^j \subset \mathcal{D}^{i+j}~~ \text{~and~~}~~[\mathcal{D}^i, \mathcal{D}^j]_c \subset \mathcal{D}^{i+j-1},$$  where $\cdot$ denotes the multiplication of $\mathcal{D}$ and $[\cdot, \cdot]_c$ is the commutator bracket and by convention $\mathcal{D}^{-1}=\{0\}$.

            A morphism $\psi: (\mathcal{D}, \mathcal{D}^i) \rightarrow (\tilde{\mathcal{D}}, \tilde{\mathcal{D}}^i)$ of quantum Poisson algebroids is given as follows.
			\begin{itemize}
				\item The map $\psi$ is a filtration preserving homomorphism, i.e. $\psi(\mathcal{D}^i) \subset \tilde{D}^i$ for $i \in \mathbb{N} \cup \{0\}$,
				\item The map $\psi: \mathcal{D} \rightarrow \tilde{\mathcal{D}}$ is a $\mathbb{K}_X$-algebra homomorphism.
			\end{itemize}
            With this notion of morphism, quantum Poisson algebroids over $X$ constitute a category, denote it by $\mathcal{QP}_X$.
            \end{Def}
			Notice that, $\mathcal{D}^0$ is a commutative associative $\mathbb{K}_X$-algebra $(\mathcal{D}^0 \cdot \mathcal{D}^0 \subset \mathcal{D}^{0}$ and $[\mathcal{D}^0, \mathcal{D}^0]_c=0)$, and 
            the morphism $\psi$ is a $\mathcal{D}^0$-module homomorphism where $\mathcal{D}^0$ acts on $\tilde{D}$ via $\psi|_{\mathcal{D}^0}: \mathcal{D}^0 \rightarrow \tilde{\mathcal{D}}^0$.
            
			\subsection{Quantum Poisson algebroid structure on Universal enveloping algebroid}
			    For a Lie algebroid $(\mathcal{O}, \mathcal{L})$, the universal enveloping algebroid $\mathscr{U}(\mathcal{O},\mathcal{L})$ has a canonical quantum Poisson algebroid structure, given as follows. Recall that, $\mathcal{D}:=\mathscr{U}(\mathcal{O},\mathcal{L})$ is an associative filtered $\mathbb{K}_X$-algebra with the canonical filtration described in (\ref{natural filtration}) as $$ \mathcal{D}^0 \subset \mathcal{D}^1 \subset \cdots \subset \mathcal{D}^n \subset \cdots \subset \mathcal{D},$$ 
				where $\mathcal{D}^n:=\mathscr{U}_{n}(\mathcal{O}, \mathcal{L})$, for $n \in \mathbb{N} \cup \{0\}$. The required compatibility conditions follow from the relation (\ref{relations in U}), as the sheaf of associative algebras $\mathcal{D}$ is generated by $\mathcal{O}$ and $\iota_{\mathcal{L}} (\mathcal{L})$, where $\iota_{\mathcal{L}}: \mathcal{L} \rightarrow \mathcal{D}$ is the canonical map.
                
           If we consider $\mathcal{D}:=\mathscr{U}(\mathcal{O}, \mathcal{L})$ with the canonical filtration $(\mathcal{D}^i)_{i \in \mathbb{N} \cup \{0\}}$ for some locally free Lie algebroid $(\mathcal{O}, \mathcal{L})$, then we have the canonical PBW isomorphism (\ref{PBW})
            $$\mathscr{S}_{\mathcal{O}}\mathcal{L} \cong \oplus_i \mathcal{D}^i/{\mathcal{D}^{i-1}}=: \mathscr{S}(\mathcal{D}).$$ The associated sheaf of symmetric algebras $\mathscr{S}_{\mathcal{O}}\mathcal{L}$ has a sheaf of classical Poisson algebras structure, given by
            the canonical Poisson bracket $\{\cdot, \cdot \}$ (sheafifying the symmetric Schouten bracket) on $\mathscr{S}_{\mathcal{O}}\mathcal{L}$ as
		\begin{equation*} \label{Poisson bracket}
		    \{D_1\otimes \cdots \otimes D_n, D'_1 \otimes \cdots \otimes D'_m\}=\sum_{i,j}[D_i,D'_j] \otimes D_1 \otimes \cdots \otimes \hat{D_i} \otimes \cdots \otimes \hat{D'_j} \otimes \cdots \otimes D'_m,
		\end{equation*}
		where $D_1, \cdots, D_n,~ D'_1, \cdots, D'_m \in \mathcal{L}$. This bracket  provides a sheaf of Poisson algebras structure on $\mathscr{S}_{\mathcal{O}}\mathcal{L}$, referred as the \emph{classical Poisson algebroid} of $\mathcal{L}$. One can view  the quantum Poisson algebroid as a deformation of the classical Poisson algebroid by a filtered quantization using the PBW theorem (\ref{PBW}) (see \cite{TS}).
\begin{Rem} For the Lie algebroid $\mathcal{L}= \mathcal{T}_X$ over a smooth variety $(X, \mathcal{O}_X)$, the sheaf of polynomial functions on the cotangent bundle $T^*X$ of $X$ is 
  \label{quantum-classical tangent case}
    $\mathcal{O}_{T^*X} \cong gr(\mathcal{D}_X)\cong \mathscr{S}_{\mathcal{O}_X} \mathcal{T}_X.$ 
 It induces a canonical Poisson algebra structure on $\mathcal{O}_{T^*X}$ or a Poisson variety structure on the cotangent bundle $T^*X$ canonically (see \cite{TS}). 
     \end{Rem}
    
 \subsection{Sheaf of differential operators on an $\mathcal{O}$-module and associated quantum Poisson algebroid} \label{Differ operator}
We consider a global counterpart of the notion of algebra of differential operators on a module \cite{LM}. This also generalizes the notion of sheaf of algebras of differential operators on a commutative ring \cite{HM, TS}.

Let $\mathcal{O}$ be an commutative associative unital $\mathbb{K}_X$-algebra and $\mathcal{E}$ be an $\mathcal{O}$-module. We consider the sheaf of filtered algebras of differential operators on $\mathcal{E}$, denoted by $\mathcal{D}iff(\mathcal{E})$ and defined inductively as follows
$$\mathcal{D}iff(\mathcal{E}):= \cup_{n \geq 0} ~\mathcal{D}iff^{\leq n}(\mathcal{E}),$$
where
$\mathcal{D}iff^{0}(\mathcal{E}):= \mathcal{E}nd_{\mathcal{O}} (\mathcal{E}),$
and for $n \in \mathbb{N}$ consider the $\mathcal{O}$-modules
\begin{equation} \label{diff operators on E}
\mathcal{D}iff^{\leq n}(\mathcal{E}):= \{\tilde{D} \in \mathcal{E}nd_{\mathbb{K}_X}(\mathcal{E})~|~ [\tilde{D}, f]_c \in \mathcal{D}iff^{\leq n-1}(\mathcal{E}),~ f \in \mathcal{O}\}.    
\end{equation}
Note that, for any $\tilde{D} \in \mathcal{D}iff^{\leq n}(\mathcal{E})$, $[\tilde{D}, f]_c:=\tilde{D} \circ f - f \circ \tilde{D}: \mathcal{E} \rightarrow \mathcal{E}$ is given by $s \mapsto \tilde{D}(f~s)  - f ~ \tilde{D}(s)$. 
By considering the standard convention $\mathcal{D}iff^{\leq -1}(\mathcal{E})=0$, one can get $\mathcal{D}iff^{0}(\mathcal{E})= \mathcal{E}nd_{\mathcal{O}} (\mathcal{E})$.

The sheaf of differential operators $\mathcal{D}(\mathcal{E}) := \mathcal{D}iff(\mathcal{E})$ on an $\mathcal{O}$-module $\mathcal{E}$ generalizes $\mathcal{D}_X=\mathcal{D}(\mathcal{O}_X)$, where $(X, \mathcal{O}_X)$ is a smooth ringed space, and carries a canonical quantum Poisson algebroid structure as follows.

\begin{Thm}  \label{diff on E}
    For an $\mathcal{O}$-module $\mathcal{E}$, the sheaf $\mathcal{D}(\mathcal{E})$ forms a quantum Poisson algebroid over $X$. 
\end{Thm}
\begin{proof}
    Inductively, we show that the sheaf $\mathcal{D}(\mathcal{E})$ forms a quantum Poisson algebroid as follows.
We work at the level of sections of the $\mathcal{O}$-submodule  $\mathcal{D}(\mathcal{E}) \subset \mathcal{E}nd_{\mathbb{K}_X}(\mathcal{E})$ and use the associated sheaf homomorphism $\circ$.
    
We use the associativity property of the product $\circ$ in $\mathcal{E}nd_{\mathbb{K}_X}(\mathcal{E})$.
First, we consider $D, D' \in \mathcal{D}^{\leq 1}(\mathcal{E})$ and show that $D \circ D' \in \mathcal{D}^{\leq 2}(\mathcal{E})$,~ $[D, D']_c \in \mathcal{D}^{\leq 1}(\mathcal{E})$ holds.
Thus, for $D, D' \in \mathcal{D}^{\leq 1}(\mathcal{E})$, we show that $[D \circ D', f]_c \in \mathcal{D}^{\leq 1}(\mathcal{E})$ and $[[D, D']_c, f]_c \in \mathcal{D}^{0}(\mathcal{E})$ for all $f \in \mathcal{O}$.
  Given $s \in \mathcal{E}$, the following identity
				$$[D \circ D', f]_c (s)
				=D(D'(fs))-f(D(D'(s)))
				=D ([D', f]_c (s))+ [D,f]_c  (D'(s)),$$
 i.e. $[D \circ D', f]_c =D \circ [D', f]_c+ [D,f]_c \circ D'$  holds for all $f \in \mathcal{O}$. Using this identity, we get
	$$[[D, D']_c, f]_c= [D \circ D', f]_c-[D' \circ D, f]_c= [D, [D',f]_c]_c+ [[D, f]_c, D']_c,$$
(the Jacobi identity). Note that, $D \circ \phi,~ \phi \circ D \in \mathcal{D}^{\leq 1}(\mathcal{E})$ and $[D, \phi]_c \in \mathcal{D}^{0}(\mathcal{E})$ for $D \in \mathcal{D}^{\leq 1}(\mathcal{E}), ~\phi \in \mathcal{D}^0(\mathcal{E})$. These implies that $[D \circ D', f]_c \in \mathcal{D}^{\leq 1}(\mathcal{E})$ and $[[D, D']_c, f]_c \in \mathcal{D}^{0}(\mathcal{E})$ for all $f \in \mathcal{O}$ (by the above identities). 

Suppose for any $D \in \mathcal{D}^{\leq i}(\mathcal{E})$ and $D' \in \mathcal{D}^{\leq j}(\mathcal{E})$ with $i, j \geq 1$ and $2 \leq i+j \leq n$, we have 
$$D \circ D' \in \mathcal{D}^{\leq i+j}(\mathcal{E})~~\text{and}~~[D, D']_c \in \mathcal{D}^{\leq i+j-1}(\mathcal{E}).$$
Then, for $\tilde{D} \in \mathcal{D}^{\leq l}(\mathcal{E})$ and $\tilde{D}' \in \mathcal{D}^{\leq m}(\mathcal{E})$ with $l, m \geq 1, ~~2 \leq l+m \leq n+1$ and for any $f \in \mathcal{O}$,  we get the following identity
$$[\tilde{D} \circ \tilde{D}', f]_c
=\tilde{D} \circ (\tilde{D}' \circ f)- (f \circ \tilde{D}) \circ \tilde{D}'
= \tilde{D} \circ [\tilde{D}', f]_c + [\tilde{D}, f]_c \circ \tilde{D}',$$
using the associative property of $\circ$ in $\mathcal{E}nd_{\mathbb{K}_X}(\mathcal{E})$ and adding subtracting suitable terms, where $[\tilde{D}, f]_c \in \mathcal{D}^{\leq l-1}(\mathcal{E})$, $[\tilde{D}', f]_c  \in \mathcal{D}^{\leq m-1}(\mathcal{E})$. Thus, $[\tilde{D} \circ \tilde{D}', f]_c \in \mathcal{D}^{\leq l+m-1}(\mathcal{E})$ for any $f \in \mathcal{O}$, implies 
$$\tilde{D} \circ \tilde{D}'  \in \mathcal{D}^{\leq l+m}(\mathcal{E})\subseteq \mathcal{D}^{\leq n+1}(\mathcal{E}).$$ 
Using the Jacobi identity for the commutator Lie bracket $[\cdot, \cdot]_c$ on $\mathcal{E}nd_{\mathbb{K}_X}(\mathcal{E})$ we get
\begin{align*}
[[\tilde{D}, \tilde{D}']_c, f]_c
=&~[\tilde{D}, [\tilde{D}',f]_c]_c+ [[\tilde{D}, f]_c, \tilde{D}']_c.
\end{align*}
           
  Here, both $[\tilde{D}, [\tilde{D}',f]_c]_c$ and $[[\tilde{D}, f]_c, \tilde{D}']_c$ are in $\mathcal{D}^{\leq l+m-2}(\mathcal{E})$ by the assumption. Thus  
  $$[\tilde{D}, \tilde{D}']_c \in \mathcal{D}^{\leq l+m-1}(\mathcal{E}) \subseteq \mathcal{D}^{\leq n}(\mathcal{E}).$$ Hence, by the induction hypothesis the desired result holds.  \end{proof}
\begin{Rem} A canonical almost commutative filtered algebra (or, quantum Poisson algebra) structure on the algebra of differential operators $Diff(E)$ on an $R$-module $E$ is mentioned 
     in \cite[Exercise 2.1.13]{VGinz2}. This is equivalent to the result Theorem \ref{diff on E} in the associated affine scheme case.
\end{Rem}  
\begin{Exm} [Sheaf of differential operators for a logarithmic foliation] Consider the analytic subspace $(Y:=V(\mathcal{I}),~ \mathcal{O}_Y:=\mathcal{O}_X/{\mathcal{I}})$  of a complex manifold  $(X, \mathcal{O}_X)$, associated to a principal ideal sheaf $\mathcal{I} \subset \mathcal{O}_X$ (see Example \ref{Foliation}). Then we have the canonical $\mathcal{O}_X$-module structure on $\mathcal{O}_Y$ and the associated quantum Poisson algebroid is $\mathcal{D}_Y:= \mathcal{D}iff(\mathcal{O}_Y)$. Moreover, we have the sheaf of logarithmic differential operators $\mathcal{D}_X(-log~Y)$, with canonical quantum Poisson algebroid structure induced from the structure of $\mathcal{D}_X$ (see Example \ref{sheaf of log diff operators}). 
Note that, if $Y$ is a free divisor of a complex manifold $X$, then we get canonical isomorphisms
$$\mathcal{D}_Y \cong \mathscr{U}(\mathcal{O}_Y, \mathcal{T}_Y) ~\text{and}~ \mathcal{D}_X(-log~Y) \cong \mathscr{U}(\mathcal{O}_X, \mathcal{T}_X(-log ~Y)).$$ In this situation the surjective Lie algebroid homomorphism
$$\rho: (\mathcal{O}_X, \mathcal{T}_X(-log~Y)) \rightarrow (\mathcal{O}_Y, \mathcal{T}_Y)$$ canonically extends to a surjective morphism of quantum Poisson algebroids $$\tilde{\rho}: \mathcal{D}_X(-log~Y) \rightarrow \mathcal{D}_Y.$$ 
\end{Exm}

			\subsection{The Lie algebroid induced from a quantum Poisson algebroid}
			For a quantum Poisson algebroid $(\mathcal{D},\mathcal{D}^i)$, we show  there is a canonical Lie algebroid structure on $(\mathcal{D}^0, \mathcal{D}^1 _c)$,  where $\mathcal{D}^0$ with the commutative $\mathbb{K}_X$-algebra structure and $\mathcal{D}^1_c := \mathcal{D}^1/{\mathcal{D}^0}$ is the quotient $\mathcal{D}^0$-module with the induced structures from $\mathcal{D}^1$.  The Lie algebroid structure on $\mathcal{D}^1_c$ is induced from the canonical Lie algebroid structure of $\mathcal{D}^1$ as follows.

				For a quantum Poisson algebroid $(\mathcal{D},\mathcal{D}^i)$ over $X$, denote the $\mathbb{K}_X$-algebra $\mathcal{D}^0$ by $\mathcal{O}$.
                The conditions $\mathcal{O} \cdot \mathcal{D}^i\subset \mathcal{D}^{i}$ and $ \mathcal{D}^i \cdot \mathcal{O} \subset \mathcal{D}^{i}$ provide an $\mathcal{O}$-bimodule structure on each $\mathcal{D}^{i}$, for $i \in \mathbb{N} \cup \{0\}$. Moreover, the condition  $[\mathcal{D}^1, \mathcal{D}^1]_c \subset \mathcal{D}^{1}$ provides a $\mathbb{K}_X$-Lie algebra structure on $ \mathcal{D}^{1}$. Thus, $\mathcal{D}^1$ forms an $\mathcal{O}$-module and $\mathbb{K}_X$-Lie algebra, induced from the quantum Poisson algebroid structure on $\mathcal{D}=\cup ^{\infty}_{i=0}\mathcal{D}^i$. This two algebraic structure are compatible by the following Leibniz rule:
				$$[D, f \cdot D']_c= f \cdot [D,D']_c+ [D, f]_c \cdot D'.$$
				Notice that $[D, -]\in \mathcal{D}er_{\mathbb{K}_X}(\mathcal{O})$ for $D \in \mathcal{D}^1$. Define the map 
				$$\mathfrak{a'}:\mathcal{D}^1 \rightarrow \mathcal{D}er_{\mathbb{K}_X}(\mathcal{O}),$$
				given by $\mathfrak{a'}(D)(f):=[D,f]_c$ for $D \in \mathcal{D}^1, ~ f \in \mathcal{O}_X$. Hence, $(\mathcal{D}^1,[\cdot, \cdot]_c, \mathfrak{a'})$ is a Lie algebroid over $(X, \mathcal{O})$.

				\begin{Thm} \label{quant-Lie} Let $(\mathcal{D}, \mathcal{D}^i)$ be a quantum Poisson algebroid.
				Then  $(\mathcal{D}^0, \mathcal{D}^1 _c)$ forms a Lie algebroid.
			\end{Thm}
			\begin{proof}
			Consider the pair $(\mathcal{D}^0, \mathcal{D}^1 _c)$ for a given quantum Poisson algebroid $(\mathcal{D}, \mathcal{D}^i)$. Thus, $\mathcal{D}^1_c$ canonically induces a $\mathbb{K}_X$-Lie algebra and an $\mathcal{O}$-module $($where $\mathcal{O}:= \mathcal{D}^0)$ structure from the algebraic structures of $\mathcal{D}^1$ as follows:
			$$[D+ \mathcal{O}, D'+\mathcal{O}]:=[D, D']_c+ \mathcal{O}, ~ f~(D+\mathcal{O}):=f\cdot D+ \mathcal{O},$$ for $D, D' \in \mathcal{D}^1, ~f \in \mathcal{O}$.
				Hence, the canonical surjective map 
                $$\rho: \mathcal{D}^1 \rightarrow \mathcal{D}^1_c ~~\text{~~defined by}~~\rho(D)=\bar{D}$$ where $\bar{D}:=D + \mathcal{O}$,
                is an $\mathcal{O}$-linear map as well as $\mathbb{K}_X$-Lie algebra homomorphism.

			Consider the $\mathcal{O}$-module homomorphism
			\begin{center}
				$\mathfrak{a}: \mathcal{D}^1_c \rightarrow \mathcal{D}er_{\mathbb{K}_X}(\mathcal{O})$ defined by $\bar{D} \mapsto \tilde{D}$ 
			\end{center}such that $\tilde{D}(f):=[D,f]_c \in \mathcal{O}$ for any $f \in \mathcal{O}$, where $\rho(D)=\bar{D}$ for some $D \in \mathcal{D}^1$.    To check the well definedness of the above assignment, consider two element $D, D' \in \mathcal{D}^1$ such that $\rho(D)= \bar{D} =\rho(D')$. Thus, $\rho (D-D')=0$ and it implies $D-D' \in \mathcal{O}$. Hence, $[D-D', f]_c=0$ for any $f \in \mathcal{O}$, implies $[D,f]_c=\mathfrak{a}(\bar{D})(f)=[D',f]_c$.  The map $\mathfrak{a}$ forms anchor for the Lie algebroid structure on $\mathcal{D}^1_c $ as follows.

			Let $\bar{D}, \bar{D'} \in \mathcal{D}^1_c$  and $f\in \mathcal{O}$. Then the following identities hold
$$[\tilde{D},\tilde{D'}]_c(f) = \tilde{D}([D',f]_c) - \tilde{D'}([D,f]_c) = [D,[D',f]_c]_c - [D',[D,f]_c]_c = [[D,D']_c,f]_c = \widetilde{[\bar{D},\bar{D'}]}(f),$$
and $$[\bar{D},f ~ \bar{D'}]
=\overline{[D, f \cdot D']_c}
=\overline{f \cdot [D, D']_c+ \mathfrak{a}'(D)(f) D'}
=f~[\bar{D}, \bar{D'}]+\mathfrak{a}(\bar{D})(f) \bar{D'}.$$
In the underlying process of the above constructions, one works with spaces of sections over open sets of $X$, along with the compatibility properties of restriction maps.
\end{proof}
Given an $\mathcal{O}$-module $\mathcal{E}$, from the above discussion we can deduce that the pair $(\mathcal{E}nd_{\mathcal{O}}(\mathcal{E}), \mathcal{D}iff^{\leq 1}(\mathcal{E}))$ has a canonical Lie algebroid structure induced from the quantum Poisson algebroid structure on $\mathcal{D}iff(\mathcal{E})$ $($see (\ref{diff operators on E})$)$. Here, we show that the pair $(\mathcal{O}, \mathcal{D}iff^{\leq 1}(\mathcal{E}))$  induces a canonical Lie algebroid structure as follows.
\begin{Prop}
   The $\mathcal{O}$-module  $\mathcal{D}iff^{\leq 1}(\mathcal{E})$ forms a Lie algebroid over $X$, given by the Lie bracket $[\cdot, \cdot]_c$ and the canonical $($anchor$)$ map
\begin{center}
  $\tilde{\sigma}: \mathcal{D}iff^{\leq 1}(\mathcal{E}) \rightarrow \mathcal{D}er_{\mathbb{K}_X}(\mathcal{O}) $\\
  ${D} \mapsto p_1([{D}, -]_c)$,
\end{center} 
where $p_1: \mathcal{D}er_{\mathbb{K}_X}(\mathcal{O}) \otimes_{\mathcal{O}} \mathcal{E}nd_{\mathcal{O}} (\mathcal{E}) \rightarrow \mathcal{D}er_{\mathbb{K}_X}(\mathcal{O})$ is the  canonical map $($uniquely$)$ determined by the universal property of the tensor product.
\end{Prop}
\begin{proof} First, we show that there is a map (forms an $\mathcal{O}$-linear and $\mathbb{K}_X$-Lie algebra homomorphism)
$$\sigma: \mathcal{D}iff^{\leq 1}(\mathcal{E}) \rightarrow \mathcal{D}er_{\mathbb{K}_X}(\mathcal{O})\otimes_{\mathcal{O}} \mathcal{E}nd_{\mathcal{O}} (\mathcal{E})$$ defined by $\sigma(D)(f \otimes s)=\sigma_D(f)(s)$, for any ${D} \in \mathcal{D}iff^{\leq 1}(\mathcal{E})$ and $f \in \mathcal{O}$, where $\sigma_D(f):=[{D}, f]_c \in \mathcal{E}nd_{\mathcal{O}} (\mathcal{E})$.
Given a section ${D} \in \mathcal{D}iff^{\leq 1}(\mathcal{E})$, it satisfies the following Leibniz rule:
$$D(f~s)=f~D(s) + \sigma_D(f) (s),$$
for $f \in \mathcal{O}$,~ $s \in \mathcal{E}$.
In the following, we show that $\sigma_D:= [D,-]_c \in \mathcal{D}er_{\mathbb{K}_X}(\mathcal{O})\otimes_{\mathcal{O}} \mathcal{E}nd_{\mathcal{O}} (\mathcal{E})$, for $D \in \mathcal{D}iff^{\leq 1}(\mathcal{E})$.

  Let $D \in \mathcal{D}iff^{\leq 1}(\mathcal{E})$,~ $f, g \in \mathcal{O}$, and $s \in \mathcal{E}$. 
 Then, the following equalities hold
 \begin{align*}
				D(f~gs)- fg~ D(s) 
				=&~f~D(gs)+ \sigma_D(f)~gs- fg~D(s)\\
				=&~ f(g~D(s) + \sigma_D (g)~s) + \sigma_D(f)~gs- fg~D(s)\\
				=&~(\sigma_D(f)~g + f~ \sigma_D (g)) (s).
			\end{align*}
This implies that, $[D, fg]_c(s)= (\sigma_D(f)~g + f~ \sigma_D (g)) (s)$ holds, and thus
$\sigma_D(fg)= \sigma_D(f)~g + f~ \sigma_D (g)$. Thus, we have
$\sigma(D)(- \otimes s) \in \mathcal{D}er_{\mathbb{K}_X}(\mathcal{O})$ and $\sigma(D)(f \otimes -) \in \mathcal{E}nd_{\mathcal{O}}(\mathcal{E})$, and hence the map $\sigma$ exists.

By definition, we obtain the desired compatibility condition (i.e. the Leibniz rule) for $\mathcal{D}iff^{\leq 1}(\mathcal{E})$, and the map $\sigma$
forms a Lie algebroid homomorphism, where $\mathcal{E}nd_{\mathcal{O}} (\mathcal{E})$ is regarded as a trivial Lie algebroid. Hence, the map $\tilde{\sigma}=p_1 \circ \sigma$ is also a Lie algebroid homomorphism.
\end{proof}
Therefore, the Lie algebroid $\mathcal{D}iff^{\leq 1}(\mathcal{E})$ is naturally related to the Atiyah algebroid $\mathcal{A}t(\mathcal{E})$.
    \subsection{Categorical correspondence between $\mathcal{LA}_X$ and $\mathcal{QP}_X$.}
		 Consider the categories $\mathcal{LA}_X$ of Lie algebroids (Definition \ref{Category LA}) and $\mathcal{QP}_X$ of quantum Poisson algebroids (Definition \ref{Category QP}) over a topological space $X$. We show
         a categorical correspondence between $\mathcal{LA}_X$ and $\mathcal{QP}_X$.
			\begin{Thm} \label{LA-QP sheaf}
			 The following assignments forms adjoint functors between the categories $\mathcal{LA}_X$ and $\mathcal{QP}_X$.
			\begin{center}
				$\mathscr{U}:\mathcal{LA}_X \rightarrow \mathcal{QP}_X$ given by $(\mathcal{O}, \mathcal{L}) \mapsto \mathscr{U}(\mathcal{O}, \mathcal{L}),$
               \end{center}
               \hspace{.3 cm}
               \begin{center}
				$\mathscr{D}^1:\mathcal{QP}_X\rightarrow \mathcal{LA}_X$ given by $(\mathcal{D},\mathcal{D}^i) \mapsto (\mathcal{D}^0, \mathcal{D}^1_c)$.
			\end{center}
			\end{Thm}
			\begin{proof}
			For an object $(\mathcal{O}, \mathcal{L})$ in the category $\mathcal{LA}_X$, we have the associated object $\mathscr{U}(\mathcal{O}, \mathcal{L})$ in $\mathcal{QP}_X$. On the other hand, for an object $(\mathcal{D}, \mathcal{D}^i)$ in $\mathcal{QP}_X$ we have the associated object $(\mathcal{D}^0, \mathcal{D}^1/{\mathcal{D}^0})$ in $\mathcal{LA}_X$.
			
			Let $(\phi_0,\phi_1): (\mathcal{O}, \mathcal{L}) \rightarrow (\mathcal{O}', \mathcal{L}')$ be a morphism in the category $\mathcal{LA}_X$. Using the properties of the universal enveloping algebroid of a Lie algebroid, we get a canonical map
			$$ \mathscr{U}(\phi_0,\phi_1)=\tilde{\phi}: \mathscr{U}(\mathcal{O}, \mathcal{L}) \rightarrow \mathscr{U}(\mathcal{O}', \mathcal{L}') $$
			defined by $$\tilde{\phi}(f ~\overline{D_1} \cdots \overline{D_n})= \phi_0(f) ~\overline{\phi_1(D_1)} \cdots \overline{\phi_1(D_n)}$$ for $f \in \mathcal{O}$, $D_1, \dots, D_n \in \mathcal{L}$, where $\overline{D'_i}=\iota_{\mathcal{L}''}(D'_i)$ given by the canonical map 
   $\iota_{\mathcal{L}''}: \mathcal{L}'' \rightarrow \mathscr{U}(\mathcal{O}'', \mathcal{L}'')$ for a Lie algebroid  $(\mathcal{O}'', \mathcal{L}'')$. We show this map is a homomorphism of sheaves of $\mathbb{K}$-algebras. Since from the construction (\ref{Universal enveloping algebroid}) the map $\tilde{\phi}(U)$ on the space of sections over an open set $U$ is a $\mathbb{K}$-algebra homomorphism, it is sufficient to check the compatibility condition, i.e. for each open set $V \subset U$  the following holds
   \begin{align} \label{compatibility on U(L)}
     res_{UV}^{\mathscr{U} (\mathcal{O}', \mathcal{L}')} \circ \tilde{\phi}(U)= \tilde{\phi}(V) \circ res_{UV}^{\mathscr{U} (\mathcal{O}, \mathcal{L})}.  
   \end{align}
   Since the map $(\phi_0, \phi_1)$ is a morphism in $\mathcal{LA}_X$, thus for $f \in \mathcal{O}(U)$, $D \in \mathcal{L}(U)$ we have
   $$res_{UV}^{\mathcal{L}'} (\phi_1(U)(f~D))= (\phi_0(U)(f))|_V~ (\phi_1(U)(D))|_V= \phi_0(V)(f|_V)~\phi_1(V)(D|_V)= \phi_1(V)(res_{UV}^{\mathcal{L}}(f~D)).$$
   Applying this identity among the generators, the required identity (\ref{compatibility on U(L)}) holds.
   Moreover, the map $\tilde{\phi}$ is a filtration preserving map, and hence a morphism  in the category $\mathcal{QP}_X$.
			
			Let $\phi: (\mathcal{D}, \mathcal{D}^i) \rightarrow (\tilde{\mathcal{D}}, \tilde{\mathcal{D}}^i)$ be a morphism in $\mathcal{QP}_X$ and $\phi|_{\mathcal{D}^i}:=\phi_i$ for $i \geq 0$. Then we show that the pair of maps $$\mathscr{D}^1(\phi)=(\phi_0,{\phi}^c_1): (\mathcal{D}^0, \mathcal{D}^1_{c}) \rightarrow (\tilde{\mathcal{D}}^0, \tilde{\mathcal{D}}^1_{c})$$ where 
			$\mathcal{D}^1_{c}=\mathcal{D}^1/{\mathcal{D}^0}$ and $\tilde{\mathcal{D}}^1_{c}=\tilde{\mathcal{D}}^1/{\tilde{\mathcal{D}}^0}$ and ${\phi}^c_1(\bar{D})=\overline{\phi_1(D)}$ for $\bar{D}  \in  \mathcal{D}^1_{c}$ associated to $D \in \mathcal{D}^1$, provides a morphism in the category $\mathcal{LA}_X$. 
            Let $\mathfrak{a}$, $\mathfrak{a}'$ are the associated anchor maps respectively. Notice that the restriction map $\phi_1$ and the induced quotient map $\phi^c_1$ are sheaf homomorphisms.
			Since the map $\phi$ is a filtration preserving $\mathbb{K}_X$-algebra homomorphism, thus we get the following identities: 
			$$\phi^c_1(f~\bar{D})=\overline{\phi_1(f \cdot {D})}=\phi_0(f)~ \overline{\phi_1(D)}=\phi_0(f)~\phi^c_1(\bar{D}),$$
			$$\mathfrak{a}'(\phi^c_1(\bar{D}))(\phi_0(f))=[\phi_1(D), \phi_0(f)]_c=\phi_0([D,f]_c)=\phi_0(\mathfrak{a}(\bar{D})(f)),$$
			for $f \in \mathcal{O}$ and $\bar{D} \in \mathcal{D}^1_{c}$. Since for $D, D' \in \mathcal{D}^1$, we have
            $$\phi_1([D, D']_c)=\phi(D~D'-D'~D)=\phi_1(D)~\phi_1(D')-\phi_1(D')~\phi_1(D)=[\phi_1(D),\phi_1(D')]_c,$$ 
			i.e. $\phi_1: \mathcal{D}^1 \rightarrow \tilde{\mathcal{D}}^1$ is a Lie algebra homomorphism, induces the canonical Lie algebra homorphism $$\phi^c_1:\mathcal{D}^1_{c} \rightarrow \tilde{\mathcal{D}}^1_{c}.$$ 
			From the above ideas, it follows that both $\mathscr{U}$ and $\mathscr{D}^1$ forms functors between the corresponding categories.

			\emph{On adjunction, unit and counit.}
			To show adjunction for the above described functors we need to show that for objects $(\mathcal{O}, \mathcal{L})$ and $(\mathcal{D}, \mathcal{D}^i)$ in  $\mathcal{LA}_X$ and $\mathcal{QP}_X$ respectively, there is a bijection between the respective morphism sets
			$$\mathscr{H}om_{\mathcal{LA}_X}((\mathcal{O}, \mathcal{L}),~(\mathcal{D}^0,\mathcal{D}^1/{\mathcal{D}^0})), ~ \mathscr{H}om_{\mathcal{QP}_X}(\mathscr{U}(\mathcal{O}, \mathcal{L}),\mathcal{D}),$$
			 such that the family of bijections is natural in $(\mathcal{O}, \mathcal{L})$ and $(\mathcal{D}, \mathcal{D}^i)$. 
			
			Using the universal property of universal enveloping algebra (see Remark \ref{Uni prop}) for given a morphism between $(\mathcal{O}, \mathcal{L})$  and $(\mathcal{D}^0,\mathcal{D}^1/{\mathcal{D}^0})$ we  get a canonical morphism between $\mathscr{U}(\mathcal{O}, \mathcal{L})$ and $\mathcal{D}$ in the respective categories. On the other hand given a morphism from $\mathscr{U}(\mathcal{O}, \mathcal{L})$ to $\mathcal{D}$ we get a morphism from $(\mathcal{O},  ~ \mathscr{U}^1(\mathcal{O}, \mathcal{L})/\mathcal{O})$ to $(\mathcal{D}^0,\mathcal{D}^1/{\mathcal{D}^0})$ and thus the canonical map $\iota: \mathcal{L} \rightarrow \mathscr{U}^1(\mathcal{O}, \mathcal{L})\subset \mathscr{U}(\mathcal{O}, \mathcal{L})$, provides a morphism from $(\mathcal{O}, \mathcal{L})$  to $(\mathcal{D}^0,\mathcal{D}^1/{\mathcal{D}^0})$ in the respective categories. Moreover, the associated naturality condition holds.

   Equivalently, we have two natural transformations, namely the unit and counit of the adjunction as follows.
			The unit of the adjunction is
			$$\eta:Id_{\mathcal{LA}} \rightarrow \mathscr{D}^1 \circ \mathscr{U}$$
			given by the canonical morphism $\eta_{(\mathcal{O}, \mathcal{L})}: (\mathcal{O}, \mathcal{L}) \rightarrow (\mathcal{O},~\mathscr{U}^1(\mathcal{O}, \mathcal{L})/\mathcal{O})$ for an object $(\mathcal{O}, \mathcal{L})$ in  $\mathcal{LA}_X$.
  The co-unit of the adjunction is 
			$$\epsilon: \mathscr{U} \circ \mathscr{D}^1 \rightarrow Id_{\mathcal{QP}}$$
			given by the canonical map $\epsilon_{\mathcal{D}}: \mathscr{U}(\mathcal{D}^0, \mathcal{D}^1/{\mathcal{D}^0}) \rightarrow \mathcal{D}$ for an object $(\mathcal{D}, \mathcal{D}^i)$ in  $\mathcal{QP}_X$.
		\end{proof}
        \subsubsection{\textbf{Application to particular cases}} We describe Theorem \ref{LA-QP sheaf} in the smooth, holomorphic and algebraic contexts.
\medskip
        
      \noindent $(1)$ \textbf{Lie–Rinehart algebras and quantum Poisson algebras.}
Let \( L \) be a Lie–Rinehart algebra over a commutative  ring $($or, $\mathbb{K}$-algebra$)$ \( R \). Consider the corresponding affine scheme $(X, \mathcal{O}_X)$, where \( X := \operatorname{Spec}(R) \) with its structure sheaf \( \mathcal{O}_X \). By localization procedure, \( L \) gives rise to a Lie algebroid \( \mathcal{L} \) over \( (X, \mathcal{O}_X) \) \cite{CW}. Then, by Theorem~\ref{LA-QP sheaf}, we obtain a quantum Poisson algebroid \( \mathscr{U}(\mathcal{O}_X, \mathcal{L}) \). The space of global sections of this sheaf coincides with the universal enveloping algebra \( \mathcal{U}(R, L) \) of the Lie–Rinehart algebra \( L \).

Conversely, let \( (D, D^i) \) be a quantum Poisson algebra. Consider the affine scheme $(X, \mathcal{O}_X)$, where \( X := \operatorname{Spec}(D^0) \) with its sheaf of regular functions \( \mathcal{O}_X \). Localization of each $D^0$-module \( D^i \) yields sheaves \( \mathcal{D}^i \), and the collection \( \mathcal{D} := \bigcup_i \mathcal{D}^i \) forms a  quantum Poisson algebroid over \( X \). Again using Theorem~\ref{LA-QP sheaf}, one obtains a Lie algebroid \( \mathcal{D}^1_c \) over \( (X, \mathcal{O}_X) \), with global sections 
$\mathcal{D}^1_c(X) \cong \mathcal{D}^1(X)/\mathcal{D}^0(X) \cong D^1/D^0$ (see \cite{GuTw}),
providing the associated Lie–Rinehart algebra \( (D^0, D^1/{D^0}) \).
\begin{Cor}
    The above correspondence defines a one-to-one association between the objects in the category of Lie–Rinehart algebras and those in the category of quantum Poisson algebras. 
\end{Cor}
In fact, following the proof of Theorem \ref{LA-QP sheaf}, it follows that this correspondence extends to a categorical adjunction between the categories of Lie–Rinehart algebras and quantum Poisson algebras.

Let \( X = \operatorname{Spec}(R) \). The category \( \mathcal{LA}_{(X, \mathcal{O}_X)} \) of quasicoherent Lie algebroids \( (\mathcal{O}_X, \mathcal{L}) \) over \( X \) forms a subcategory of \( \mathcal{LA}_X \), and is equivalent to the category \( LR_R \) of \( (\mathbb{K}, R) \)-Lie–Rinehart algebras (see Example~\ref{LR on affine schemes}). Similarly, the category \( \mathcal{QP}_{(X, \mathcal{O}_X)} \) of quantum Poisson algebroids \( (\mathcal{D}, \mathcal{D}^i) \) over \( X \), with \( \mathcal{D}^0 = \mathcal{O}_X \) and \( \mathcal{D}^1 \) a quasicoherent \( \mathcal{O}_X \)-module, corresponds to the category \( QP_R \) of quantum Poisson algebras with degree-zero part \( R \), in the sense of \cite{GP}.
As mentioned in Remark \ref{usual morphisms}, morphisms in the categories \( \mathcal{LA}_{(X, \mathcal{O}_X)} \) and \( \mathcal{QP}_{(X, \mathcal{O}_X)} \) are taken such as the identity map is the only map on the underlying morphism from $\mathcal{O}_X$ to $\mathcal{O}_X$.

\begin{Prop}
We have the following equivalence of categories:
\[
LR_R \cong \mathcal{LA}_{(X, \mathcal{O}_X)}, \quad QP_R \cong \mathcal{QP}_{(X, \mathcal{O}_X)}.
\]
These equivalences are induced by the canonical equivalence between the category of \( R \)-modules and the category of quasicoherent \( \mathcal{O}_X \)-modules, with adjoint functors given by localization and global sections (\cite{CW}).
\end{Prop}
\begin{Cor}
 We have canonical adjoint functors between the categories $\mathcal{LA}_{(X, \mathcal{O}_X)}$ and $\mathcal{QP}_{(X, \mathcal{O}_X)}$.   
\end{Cor}
\medskip   
      \noindent $(2)$ \textbf{Smooth Lie algebroids and locally free quantum Poisson algebroids.}\\
Let \( L \) be a smooth Lie algebroid over a smooth manifold \( X \), and let \( \mathcal{L} := \Gamma_X(L) \) denote its sheaf of sections. Then \( \mathcal{L} \) is a locally free Lie algebroid of finite rank over the  ringed space \( (X, C^\infty_X) \). By an analogue of the Serre–Swan theorem for smooth vector bundles (see \cite{JNes}), the space of global sections \( \bar{L} := \mathcal{L}(X) \) is a finitely generated projective \( R := C^\infty(X) \)-module, algebraically identified with \( L \). 
The quantum Poisson algebroid \( \mathscr{U}(C^\infty_X, \mathcal{L}) \) associated to \( \mathcal{L} \) satisfies
\[
\mathscr{U}(C^\infty_X, \mathcal{L})(X) \cong \mathcal{U}(R, \bar{L}),
\]
the universal enveloping algebra of the Lie–Rinehart algebra \( (R, L) \), equipped with its canonical filtration. In this setting, the filtered pieces satisfy \( \mathcal{D}^0(X) \cong R \), \( \mathcal{D}^1(X) \cong R \oplus L \), and thus  \( \mathcal{D}^1_c(X) := \mathcal{D}^1(X)/\mathcal{D}^0(X) \cong L \).

Conversely, let \( (D, D^i) \) be a quantum Poisson algebra with \( D^0 = C^\infty(X) \) and \( D^1 \) a finitely generated projective \( D^0 \)-module. Localizing each \( D^i \) yields a sheaf \( \mathcal{D}^i \), and \( (\mathcal{D}, \mathcal{D}^i) \) becomes a quantum Poisson algebroid over \( X \), with \( \mathcal{D}^0 = C^\infty_X \) and \( \mathcal{D}^1 \) locally free of finite rank, called a locally free quantum Poisson algebroid of finite rank. The short exact sequence
$\mathcal{D}^0 \hookrightarrow \mathcal{D}^1 \twoheadrightarrow \mathcal{D}^1_c$
then yields an exact sequence of projective modules:
\[
\mathcal{D}^0(X) \hookrightarrow \mathcal{D}^1(X) \twoheadrightarrow \mathcal{D}^1_c(X),
\]
which defines a transitive extension of Lie–Rinehart algebras. This splits as \( \mathcal{D}^1(X) \cong \mathcal{D}^0(X) \oplus \mathcal{D}^1_c(X) \), and the pair \( (\mathcal{D}^0(X), \mathcal{D}^1_c(X)) \) corresponds to the Lie algebroid \( (\mathcal{D}^0, \mathcal{D}^1_c) \).

Following the terminology mentioned in Remark \ref{usual morphisms}, we get the following result.
\begin{Cor}
Let $X$ be a smooth manifold.
 Denote the category of smooth Lie algebroids and of locally free quantum Poisson algebroids of finite rank over $(X, C^\infty_X )$ by \( \mathcal{LA}_{(X, C^\infty_X)} \)  and  \( \mathcal{QP}_{(X, C^\infty_X)} \) respectively. Then there is a induced one-to-one correspondence between these subcategories of $\mathcal{LA}_X$ and $\mathcal{QP}_X$ respectively.
\end{Cor}

\medskip

\noindent $(3)$ \textbf{Holomorphic Lie algebroids and locally free quantum Poisson algebroids.} Let \( L \) be a holomorphic Lie algebroid over a complex manifold \( X \). Consider the $\mathcal{O}_X$-module of its sheaf of sections  \( \mathcal{L} := \Gamma_X(L) \). Thus, $\mathcal{L}$ forms a locally free Lie algebroid  of finite rank over \( (X, \mathcal{O}_X) \). The associated quantum Poisson algebroid \( \mathcal{D} := \mathscr{U}(\mathcal{O}_X, \mathcal{L}) \) has \( \mathcal{D}^0 = \mathcal{O}_X \) and \( \mathcal{D}^1 \) a locally free \( \mathcal{O}_X \)-module of finite rank.

Conversely, given a quantum Poisson algebroid \( (\mathcal{D}, \mathcal{D}^i) \) over the complex manifold \( X \) with \( \mathcal{D}^1 \) is locally free of finite rank over \( \mathcal{D}^0=\mathcal{O}_X \), the quotient \( \mathcal{L} := \mathcal{D}^1/\mathcal{D}^0 \) forms a locally free Lie algebroid of finite rank over \( (X, \mathcal{O}_X) \). The fiber of the corresponding holomorphic Lie algebroid of $\mathcal{L}$ at a point \( x \in X \) is the complex vector space \( \mathcal{L}_x / \mathcal{M}_x \mathcal{L}_x \), where \( \mathcal{M}_x \) is unique the maximal ideal of the stalk \( \mathcal{O}_{X,x} \). The local triviality  of the underlying vector bundle follows from the assumption that $\mathcal{L}$ is a locally free  $\mathcal{O}_X$-module of finite rank.

Following the terminology mentioned in Remark \ref{usual morphisms}, we get the following result.
\begin{Cor}
Let $X$ be a complex manifold.
Denote the categories of holomorphic Lie algebroids and of locally free quantum Poisson algebroids of finite rank over $(X, \mathcal{O}_X)$ by \( \mathcal{LA}_{(X, \mathcal{O}_X)} \) and \( \mathcal{QP}_{(X, \mathcal{O}_X)} \), respectively. Then there is the induced adjunction between these subcategories of $\mathcal{LA}_X$ and $\mathcal{QP}_X$ respectively.
\end{Cor}

\section{Sheaves of almost commutative algebras} \label{Sec 4}
In this section, we discuss a sheaf-theoretic analogue of almost commutative filtered algebras to unify and generalize several existing frameworks, for instance, Maakestad’s almost commutative filtered (rings) algebras for Lie–Rinehart algebras \cite{HM} and Tortella’s almost polynomial filtered algebras for holomorphic Lie algebroids \cite{PT}. We observe that the sheaves of almost commutative filtered algebras are equivalent to quantum Poisson algebroids.

We recall that in \cite{HM}, Maakestad developed twisted universal enveloping algebras within the theory of filtered algebras and studied sheaves of rings of differential operators over affine schemes. In a parallel development, Tortella \cite{PT} introduced sheaves of almost polynomial filtered algebras in the holomorphic Lie algebroid setting.
Our framework includes both approaches and provides a unified perspective with new geometric and algebraic consequences. As a concrete example, we construct twisted universal enveloping algebroids of Lie algebroids over commutative ringed spaces.

\begin{Def}
   Let  $(\mathscr{U}, \mathscr{U}_i)$  be a sheaf of filtered associative $\mathbb{K}_X$-algebras over a topological space  $X$. It is said to be a sheaf of almost commutative algebras if the following conditions hold.

 \begin{enumerate}
\item The constant sheaf $\mathbb{K}_X$ is in the center of $\mathscr{U}$.

\item  $\mathscr{U}_i~ \mathscr{U}_j \subset \mathscr{U}_{i+j}$ for $i,j\in \mathbb{N} \cup \{0\}$.

\item  For $D_n \in \mathscr{U}_{i_n}$, $n=1, \dots, k$ $(k \in \mathbb{N})$, the following identity holds 
\begin{align} \label{cond 3}
    D_1 \cdots D_k - D_{\sigma(1)} \cdots D_{\sigma(k)} \in \mathscr{U}_{l-1},
\end{align}
     where $\sigma \in S_k$ $($the permutation group of $k$ elements$)$ and $l:=i_1+ \cdots +i_k$. 
\end{enumerate}
\end{Def}
We say $(\mathscr{U}, \mathscr{U}_i)$ is a sheaf of almost commutative $\mathbb{K}_X$-algebras over a ringed space $(X, \mathcal{O})$ where $\mathcal{O}:= \mathscr{U}_0$.
\begin{Rem}
    The condition (\ref{cond 3}) implies that $[\mathscr{U}_i, \mathscr{U}_j]_c \subset \mathscr{U}_{i+j-1}$ holds, for $i,j \in \mathbb{N} \cup \{0\}$. Conversely, this inclusion allows us to recover condition (\ref{cond 3}), since any permutation can be expressed as a product of transpositions. Thus, this new notion is equivalent to the notion of quantum Poisson algebroids.
\end{Rem}

\noindent \label{Simpson's axioms}
 A  sheaf of almost commutative algebras $(\mathscr{U}, \mathscr{U}_i)$ is called a \emph{Simpson's sheaf of ring of differential operators} (see \cite{PT, HM}), if in addition, it satisfies the following conditions.
 \begin{enumerate}
 \item  $(X, \mathcal{O}_X)$ is a smooth complex algebraic variety or a complex manifold and $\mathscr{U}_0=\mathcal{O}_X$,
    \item The graded $\mathcal{O}_X$-modules $\mathscr{U}_i/{\mathscr{U}_{i-1}}$ are coherent for $i \geq 1$,

\item The canonical map $\gamma: \mathscr{S}_{\mathcal{O}_X}(\mathscr{U}_1/\mathscr{U}_0) \rightarrow gr (\mathscr{U})$ is an isomorphism. 
 \end{enumerate}
In particular, in the classical case, i.e. when $\mathscr{U}_1 / \mathscr{U}_0$ is the (logarithmic) tangent sheaf $\mathcal{T}_X$ or $\mathcal{T}_X(-\log~ Y)$ for a (normal crossing) divisor $Y$ of $X$, the study was carried out by Simpson \cite{Simpson}.

Next, we construct a twisted universal enveloping algebroid of a Lie algebroid over a ringed space.
\subsection{Twisted universal enveloping algebroid} \label{Twisted enveloping algebroid}
We construct a twisted version of the universal enveloping algebroid associated with a Lie algebroid $(\mathcal{O}, \mathcal{L})$ equipped with a 2-cocycle $\omega \in \mathcal{Z}^2(\mathcal{L}, \mathcal{O})$ in the (generalized) Chevalley–Eilenberg–de Rham complex, denoted by $\mathscr{U}(\mathcal{O}, \mathcal{L}, \omega)$. This construction generalizes the framework of twisted universal enveloping algebras for Lie-Rinehart algebras described in \cite{HM}, as well as the holomorphic Lie algebroid setting in \cite{PT}. When $\omega = 0$, we recover the classical universal enveloping algebroid (see Section \ref{Universal enveloping algebroid}). Moreover, under the assumption that $\mathcal{L}$ is locally free as a left $\mathcal{O}$-module, we establish a (twisted) Poincaré–Birkhoff–Witt (PBW) theorem for $\mathscr{U}(\mathcal{O}, \mathcal{L}, \omega)$.

Let $(\mathcal{L}, [\cdot, \cdot], \mathfrak{a})$ over $(X, \mathcal{O})$ be a Lie algebroid and $(\mathcal{E}, \nabla)$ be an $\mathcal{L}$-module. 
For a cohomology class $\bar{\omega} \in \mathbb{H}^2(\mathcal{L}, \mathcal{E})$, we get an abelian Lie algebroid extension of $\mathcal{L}$ by $\mathcal{E}$ (up to isomorphism) and a sheaf of almost commutative algebras as follows.
First using the local description about Lie-Rinehart algebra case from \cite{HM}, for each open set $U \subset X$, we get a cohomology class $\omega_U \in H^2(\mathcal{L}(U), \mathcal{E}(U))$ for the Lie-Rinehart algebra $\mathcal{L}(U)$ with coefficient in $\mathcal{E}(U)$, which associates with new Lie-Rinehart algebra $\mathcal{L}(U)_{\omega_U}$ whose Lie bracket  and anchor are $[\cdot, \cdot]_{\omega_U}$ and $\mathfrak{a}_U$ respectively. We shifify these local descriptions and, using the hypercohomology functor, we get the following.
 Consider the $\mathcal{O}$-module $\mathcal{E}\oplus \mathcal{L}$ with the $\mathbb{K}_X$-Lie algebra structure given by
        $$[(s,D), (s',D')]_{\omega}:=(\nabla_D(s')-\nabla_{D'}(s)+\omega(D,D'),~ [D, D']),$$
        for $s, s'\in \mathcal{E},~D, D' \in \mathcal{L}$, together with the map
        $$\mathfrak{a}_{\omega}:\mathcal{E} \oplus \mathcal{L} \rightarrow \mathcal{D}er_{\mathbb{K}_X}(\mathcal{O}),~\text{defined by}~ (s, D) \mapsto \mathfrak{a}(D).$$
        Thus, $\mathcal{L}_{\omega}:=(\mathcal{E} \oplus \mathcal{L}, [\cdot, \cdot]_{\omega}, \mathfrak{a}_{\omega})$ forms a Lie algebroid over $(X, \mathcal{O})$, providing the abelian Lie algebroid extension (an equivalence class)
        \begin{align} \label{L_w from L}
           0 \rightarrow \mathcal{E} \rightarrow \mathcal{L}_{\omega} \rightarrow \mathcal{L} \rightarrow 0.  
        \end{align}

      Let $s$ be a generator for the free $\mathcal{O}$-module $\mathcal{E}=\mathcal{O}~s:=\{fs~|~ f \in \mathcal{O}\}$ and let  $\mathcal{E} \hookrightarrow \mathcal{L}_{\omega} \twoheadrightarrow \mathcal{L}$ be the extension of $\mathcal{L}$ by $\mathcal{E}$ corresponding to $\bar{\omega} \in \mathbb{H}^2(\mathcal{L}, \mathcal{E})$. In this case the standard $\mathcal{L}$-connection on $\mathcal{E}$ is given by the anchor map $\mathfrak{a}$ as follows:
      $$\nabla_D(f s):=\mathfrak{a}(D)(f)s,$$
for $D \in \mathcal{L}$ and $f \in \mathcal{O}$. This provides a global description of the assignment  (a presheaf) $U \mapsto \mathcal{L}(U)_{\omega_U}$.

Now, consider the sheaf of tensor algebras $\mathscr{T}(\mathcal{L}_{\omega}):= \oplus_{k \geq 0} \mathcal{L}_{\omega}^{\otimes^k_{\mathcal{O}}}$ of the sheaf of $\mathbb{K}$-Lie algebras $\mathcal{L}_{\omega}$, i.e. $\mathscr{T}(\mathcal{L}_{\omega})$ is the sheafification of the presheaf of $\mathbb{K}$-algebras $$U \mapsto T(\mathcal{L}_{\omega}(U)):= \mathbb{K} \oplus \mathcal{L}_{\omega}(U) \oplus \cdots \oplus \mathcal{L}_{\omega}(U)^{\otimes^n_{\mathcal{O}(U)}} \oplus \cdots,$$
and  the associated sheaf of universal enveloping algebras $\mathscr{U}(\mathcal{L}_{\omega})$ for the sheaf of $\mathbb{K}$-Lie algebras $\mathcal{L}_{\omega}$ is the sheafification of the presheaf $U \mapsto  T(\mathcal{L}_{\omega}(U))/I(U)$, where $I(U)$ is the two sided ideal of $T(\mathcal{L}_{\omega}(U))$ generated by $$\{\tilde{D} \otimes \tilde{D'}- \tilde{D'} \otimes \tilde{D} -[\tilde{D}, \tilde{D'}]_{\omega}~|~ \tilde{D}, \tilde{D'} \in \mathcal{L}_{\omega}(U)\}.$$ 
Thus, we have the canonical quotient map of associative $\mathbb{K}_X$-algebras
\begin{equation*}
    \mathcal{P}: \mathscr{T}(\mathcal{L}_{\omega}) \rightarrow \mathscr{U}(\mathcal{L}_{\omega}).
\end{equation*}
   Consider the $\mathcal{O}$-module $\mathscr{U}^+:= \mathcal{P}(\mathscr{T}^1(\mathcal{L}_{\omega}))$, where $\mathscr{T}^1(\mathcal{L}_{\omega}):=\oplus_{k \geq 1} \mathcal{L}_{\omega}^{\otimes^k_{\mathcal{O}}}$. Thus, we have the canonical maps
   \begin{equation} \label{O to U}
        \mathcal{P}_{\mathcal{O}}: \mathcal{O} \rightarrow \mathscr{U}^+~\text{defined by}~\mathcal{P}_{\mathcal{O}}(f)=fs,~\text{for}~ f \in \mathcal{O}, ~\text{and}
   \end{equation}
  \begin{equation} \label{L to U}
      \mathcal{P}_{\mathcal{L}}: \mathcal{L} \rightarrow \mathscr{U}^+~ \text{defined by}~ \mathcal{P}_{\mathcal{L}}(D)= \mathcal{P}(D), ~\text{for}~ D \in \mathcal{L}.
  \end{equation}
Finally, $\mathcal{P}_{\mathcal{L}_\omega}: \mathcal{L}_\omega \rightarrow \mathscr{U}^+$ defined by $\mathcal{P}_{\mathcal{L}_\omega}(\tilde{D}):= \mathcal{P}(\tilde{D})$, for $\tilde{D} \in \mathcal{L}_\omega$. Now, consider the two sided ideal $\mathcal{J}_\omega$ in $\mathscr{U}^+$ generated by
\begin{align*}
 \{\mathcal{P}_{\mathcal{L}_\omega}(f \tilde{D})- \mathcal{P}_{\mathcal{O}}(f) \mathcal{P}_{\mathcal{L}_\omega}(\tilde{D})~|~ f \in \mathcal{O}, \tilde{D} \in \mathcal{L}_\omega\}.  
\end{align*}
Thus, the universal enveloping algebroid of $(\mathcal{O}, \mathcal{L})$ twisted by $\omega \in \mathcal{Z}^2(\mathcal{L}, \mathcal{E})$, denoted by $\mathscr{U}(\mathcal{O}, \mathcal{L}, \omega)$, is 
$$\mathscr{U}(\mathcal{O}, \mathcal{L}, \omega):= \mathscr{U}^+/{\mathcal{J}_{\omega}}.$$

Let $\mathscr{U}_i(\mathcal{O}, \mathcal{L}, {\omega}):= \mathcal{P}(\mathscr{T}_i(\mathcal{L}_{\omega}))$, where $\mathscr{T}_i(\mathcal{L}_{\omega}):= \oplus^i_{k=0} \mathcal{L}_{\omega}^{\otimes^k_{\mathcal{O}}} \subset \mathscr{T}(\mathcal{L}_{\omega})$. We get a filtration of $\mathscr{U}(\mathcal{O}, \mathcal{L}, {\omega})$
$$\mathcal{O} \subset \mathscr{U}_1(\mathcal{O}, \mathcal{L}, {\omega}) \subset \cdots \subset \mathscr{U}_i(\mathcal{O}, \mathcal{L}, {\omega}) \subset \cdots ~.$$
In particular, for the case of holomorphic Lie algebroids, an equivalent approach for constructing such sheaves of algebras $($satisfying the Simpson's axioms$)$ is described in \cite{PT}.
\begin{Rem}\label{embedding of $L$ in twisted univ alg}
 Moreover, if $\mathcal{L}$ is locally free as $\mathcal{O}$-module, the the map $\mathcal{P}_{\mathcal{L}}$ is an embedding of $\mathcal{L}$ in $\mathscr{U}(\mathcal{O}, \mathcal{L}, {\omega})$.   
\end{Rem}
\begin{Rem} \label{0-twist uni alg}
    In particular, for the trivial $2$-cocycle we get, $\mathscr{U}(\mathcal{O}, \mathcal{L}, 0)= \mathscr{U}(\mathcal{O}, \mathcal{L})$, recover the standard one.
\end{Rem}
\begin{Exm} [Universal enveloping algebroids of Atiyah algebroids] Let $(X, \mathcal{O}_X)$ be a non-singular ringed space. The Atiyah sequence for an locally free $\mathcal{O}_X$-module $\mathcal{E}$ is
\begin{align} \label{Atiyah sequence}
  0 \rightarrow \mathcal{E}nd_{\mathcal{O}_X}(\mathcal{E}) \rightarrow \mathcal{A}t(\mathcal{E}) \overset{\sigma}{\rightarrow} \mathcal{T}_X \rightarrow 0. 
\end{align}
Thus, $\mathcal{A}t(\mathcal{E})$ provides an abelian Lie algebroid extension of $\mathcal{T}_X$ by $ \mathcal{E}nd_{\mathcal{O}_X}(\mathcal{E})$.
    For any $x \in X$, there is an open set $U_x \subset X$ around $x$, on the level of space of sections of (\ref{Atiyah sequence}) on $U_x$ (i.e. locally), for the associated free (or, projective) Lie-Rinehart algebras we have (see \cite{NK-PS} for the following details)
    $$\mathcal{U}(\mathcal{O}_X(U_x), \mathcal{A}t(\mathcal{E}(U_x))) \cong U({E}nd_{\mathcal{O}_X(U_x)}(\mathcal{E}(U_x)) ~\#_{\sigma}~ \mathcal{D}iff(\mathcal{O}_X(U_x)),$$
    where $\mathcal{A}t(\mathcal{E}(U_x))$ is the Lie-Rinehart algebra of $1$st order differential operators on $\mathcal{E}(U_x)$ with scalar symbols,
    $U({E}nd_{\mathcal{O}_X(U_x)}(\mathcal{E}(U_x)))$ is the universal enveloping algebra of the $\mathcal{O}_X(U_x)$-Lie algebra ${E}nd_{\mathcal{O}_X(U_x)}(\mathcal{E}(U_x))$, $\sigma$ is a suitable Hopf $2$-cocycle and $\#_{\sigma}$ denotes the corresponding crossed product. Thus, one can extend this result to the level of universal enveloping algebroids by sheafification techniques. 

    By the correspondence (\ref{Ext-H^2}), there is a cohomology class $\bar{\omega} \in \mathbb{H}^2(\mathcal{T}_X, \mathcal{E}nd_{\mathcal{O}_X}(\mathcal{E}))$ such that we get the isomorphism $\mathcal{A}t(\mathcal{E}) \cong \mathcal{E}nd_{\mathcal{O}_X}(\mathcal{E}) \oplus \mathcal{T}_X$ of Lie algebroids, where right side Lie algebroid has a Lie bracket and anchor given by the twist of $\omega$  $($i.e., $\mathcal{A}t(\mathcal{E}) \cong (\mathcal{T}_X)_{\omega})$. 
    By a canonical quotient of $\mathscr{U}(\mathcal{O}_X, \mathcal{A}t(\mathcal{E}))$, the twisted universal enveloping algebroid $\mathscr{U}(\mathcal{O}_X, \mathcal{T}_X, \omega)$ is expressed as (from the constructions)
    $$\mathscr{U}(\mathcal{O}_X, \mathcal{A}t(\mathcal{E}))/{\mathcal{J}_{\omega}} \cong \mathscr{U}(\mathcal{O}_X, \mathcal{T}_X, \omega),$$
    where $\mathcal{J}_{\omega}$ is a two sided ideal sheaf of the universal enveloping algebroid of the Lie algebroid $(\mathcal{T}_X)_{\omega}$.
\end{Exm}
In the following, we demonstrate an application of the notion of twisted universal enveloping algebroids to the study of Lie algebroid connections on an $\mathcal{O}$-module that is not necessarily flat, via their representations. In particular, we extend the correspondence result for Lie–Rinehart algebras [Proposition 4.18, \cite{HM}] to a global setting, while also relaxing the conditions in the holomorphic case treated in [Proposition 5.3, \cite{PT}].
\begin{Thm} \label{connection with prescribed curvature}
Let $\mathcal{E}$ be a left $\mathcal{O}$-module and $w \in \mathcal{Z}^2(\mathcal{L}, \mathcal{O})$ for a locally free Lie algebroid $\mathcal{L}$ over $(X, \mathcal{O})$.
     There is a one-to-one correspondence between $\mathscr{U}(\mathcal{O}, \mathcal{L}, \omega)$-module structures on $\mathcal{E}$ and $\mathcal{L}$-connections on $\mathcal{E}$ with curvature type $\omega$.
\end{Thm}
 \begin{proof}
     Consider the Lie algebroid $\mathcal{L}_{\omega}$ associated to $\mathcal{L}$ twisted by $\omega$. Consider  a left $\mathscr{U}(\mathcal{O}, \mathcal{L}, \omega)$-module structure on $\mathcal{E}$ given by the action 
     \begin{equation} \label{U action}
      (\tilde{D}, s) \mapsto \tilde{D} ~s, ~\text{for}~ \tilde{D} \in \mathscr{U}(\mathcal{O}, \mathcal{L}, \omega)~\text{and}~s \in \mathcal{E}.   
     \end{equation}
     The above action can be understood locally via the underlying action of the presheaf $\mathcal{U}(\mathcal{O}, \mathcal{L}, \omega)$ on $\mathcal{E}$, where $\mathcal{U}(\mathcal{O}, \mathcal{L}, w)$ is the presheaf of twisted universal enveloping algebras corresponding to the sheaf of Lie-Rinehart algebras $\mathcal{L}$ and $w$ is a section of the presheaf $U \mapsto H^2(\mathcal{L}(U), \mathcal{O}(U))$ associated to $\omega$. We subsequently sheafify it using the compatibility properties with restriction morphisms. 
 Using the map $\mathcal{P}_{\mathcal{L}}$ described in (\ref{L to U}) and the  action (\ref{U action}), define the $\mathcal{L}$-connection $\nabla$ on $\mathcal{E}$ defined by $$\nabla_D(s):=\mathcal{P}_{\mathcal{L}}(D)~s,$$
      where $D \in \mathcal{L}$ and $s \in \mathcal{E}$. 
      Since the Lie algebroid $\mathcal{L}$ is locally free $\mathcal{O}$-module, the above action determined by $\mathcal{P}_{\mathcal{L}}$, is well defined (see Remark \ref{embedding of $L$ in twisted univ alg}).
      Notice that, using the Lie bracket of $\mathcal{L}_{\omega}=(\mathcal{O}~z \oplus \mathcal{L}, [\cdot, \cdot]_{\omega}, \mathfrak{a}_{\omega})$ we get
      \begin{equation*}
        [(0, D), (0, D')]_{\omega}=(\omega(D, D'), [D, D']), 
      \end{equation*}
       and on sections of the sheaf of  universal enveloping algebra $\mathscr{U}(\mathcal{L}_{\omega})$, we obtain the identity
      \begin{equation*}
        D~D'-D'~D=[D, D']+\omega(D, D'), 
      \end{equation*}
     holds  for all $D, D' \in \mathcal{L}$.
      Thus, using (\ref{O to U}) and (\ref{L to U}) we get the following identity in $\mathscr{U}(\mathcal{O}, \mathcal{L}, \omega)$ 
      \begin{equation} \label{connection with curvature type}
         \mathcal{P}_{\mathcal{L}}(D)~ \mathcal{P}_{\mathcal{L}}(D') - \mathcal{P}_{\mathcal{L}}(D')~ \mathcal{P}_{\mathcal{L}}(D)= \mathcal{P}_{\mathcal{L}}([D, D'])+ \mathcal{P}_{\mathcal{O}}(\omega(D, D')), 
      \end{equation}
      holds for all $D, D' \in \mathcal{L}$.
      To determine the curvature $R_{\nabla}$ $($see (\ref{curvature type})$)$ of the $\mathcal{L}$-connection $(\mathcal{E}, \nabla)$, consider
      $$R_{\nabla}(D, D')=[\nabla_D, \nabla_{D'}]_c-\nabla_{D, D'}= \omega(D, D'),$$
      by the definition of a connection $\nabla$ and the identity (\ref{connection with curvature type}).

      Conversely, consider an $\mathcal{L}$-connection $\nabla$ on $\mathcal{E}$ with curvature type $\omega$. Thus, the $\mathcal{O}$-linear map 
      $$\nabla: \mathcal{L} \rightarrow \mathcal{E}nd_{\mathbb{K}_X}(\mathcal{E}),$$
      extends to  a canonical $\mathbb{K}_X$-algebra homomorphism
      $\theta: \mathscr{T}^1(\mathcal{L}_{\omega}) \rightarrow \mathcal{E}nd_{\mathbb{K}_X}(\mathcal{E})$
      defined by $$\theta(\otimes_i(f_i z+ D_i)):= \textstyle\prod_i(f_i ~Id_{\mathcal{E}}+ \nabla_{D_i}),$$ for $f_i \in \mathcal{O}$ and $D_i \in \mathcal{L}$, where $Id_{\mathcal{E}}$ is the identity map on $\mathcal{E}$. This map induces a homomorphism of $\mathcal{O}$-modules and of $\mathbb{K}_X$-algebras 
      $$\mathscr{U}(\mathcal{O}, \mathcal{L}, \omega) \rightarrow \mathcal{E}nd_{\mathbb{K}_X}(\mathcal{E}),$$
      and thus induces a left $\mathscr{U}(\mathcal{O}, \mathcal{L}, \omega)$-module structure on $\mathcal{E}$.
 \end{proof}
\begin{Cor}
    When we take the cohomology class $\omega=0$, we recover the classical one-to-one correspondence between representations of a Lie algebroid $\mathcal{L}$ and modules over its universal enveloping algebroid $\mathscr{U}(\mathcal{O}, \mathcal{L})$.
\end{Cor}
\begin{Rem}
    Moreover, for any locally free Lie algebroid $\mathcal{L}$ over $(X, \mathcal{O})$, the one-to-one correspondence can be extended to a categorical equivalence between the category of $\mathcal{L}$-connections of curvature type $w$ and the category of left $\mathscr{U}(\mathcal{O}, \mathcal{L},  \omega)$-modules. For the detailed local description, see [Corollary $4.20$, \cite{HM}].
\end{Rem}
Now we state and prove a twisted version of the PBW Theorem (see (\ref{PBW}) of Section \ref{Sec 2}) as follows.

\begin{Thm}\label{twisted PBW theorem} 
Let $\mathcal{L}$ be a locally free Lie algebroid over a commutative ringed space $(X, \mathcal{O})$. There is a canonical isomorphism (known as the twisted PBW map) of graded $\mathcal{O}$-algebras
$$gr(\mathscr{U}(\mathcal{O}, \mathcal{L}, \omega)) \cong \mathscr{S}_{\mathcal{O}}\mathcal{L},$$
where $gr(\mathscr{U}(\mathcal{O}, \mathcal{L}, \omega))$ is the associated sheaf of graded algebras of $\mathscr{U}(\mathcal{O}, \mathcal{L}, \omega)$ with  the ascending filtration.
\end{Thm}
\begin{proof} 
We begin by introducing some new notions, based on the earlier ones, in order to prove the theorem.
Consider the $\mathcal{O}$-module $\tilde{\mathcal{L}}:= \mathcal{O}s' \oplus \mathcal{L}_{\omega}$ with the Lie bracket
$$[f_1s'+D'_1, f_2s'+D'_2]_{\tilde{\omega}}:= (\mathfrak{a}_{\omega}(D'_1)(f_2)-\mathfrak{a}_{\omega}(D'_2)(f_1))s'+ [D'_1, D'_2]_{\omega}.$$
Since as a left $\mathcal{O}$-module $\tilde{\mathcal{L}}=\mathcal{O}s' \oplus \mathcal{O}s \oplus \mathcal{L}$ $($see (\ref{L_w from L})$)$,
Consider the map $\tilde{\mathfrak{a}}_{\omega}: \tilde{\mathcal{L}} \rightarrow \mathcal{D}er_{\mathbb{K}_X}(\mathcal{O})$ defined as $$\tilde{\mathfrak{a}}_{\omega}(f_1s'+f_2 s+ D)= \mathfrak{a}(D)$$ where $f_1s'+f_2 s+ D \in \tilde{\mathcal{L}}$. Thus, provides a Lie algebroid $(\tilde{\mathcal{L}}, [\cdot, \cdot]_{\tilde{\omega}}, \tilde{\mathfrak{a}}_{\omega})$ over $(X, \mathcal{O})$.

Consider the canonical map $\tilde{\mathcal{P}}: \mathscr{T}^1(\tilde{\mathcal{L}}) \rightarrow  \mathscr{U}(\mathcal{O}, \mathcal{L}_{\omega})$ and its restriction $\tilde{\mathcal{P}}|_{\tilde{\mathcal{L}}}: \tilde{\mathcal{L}}\rightarrow  \mathscr{U}(\mathcal{O}, \mathcal{L}_{\omega})$. Denote 
 $$\mathscr{U}(\mathcal{O}, \mathcal{L}_{\omega}, t) := \mathscr{U}(\mathcal{O}, \mathcal{L}_{\omega}) \langle t-1 \rangle,$$ where $t:= \tilde{\mathcal{P}}|_{\tilde{\mathcal{L}}}(s)$ and $\mathscr{U}(\mathcal{O}, \mathcal{L}_{\omega}) \langle t-1 \rangle$ is the ideal generated by the section $t-1$ in $\mathscr{U}(\mathcal{O}, \mathcal{L}_{\omega})$.
 It follows that $\mathscr{U}(\mathcal{O}, \mathcal{L}_{\omega}, t)$ has a canonical ascending filtration.
 
  To prove the desired result, it is enough to show that
  there is a canonical filtered $\mathbb{K}_X$-algebra isomorphism
  $$\mathscr{U}(\mathcal{O}, \mathcal{L}_{\omega}, t) \cong \mathscr{U}(\mathcal{O}, \mathcal{L}, \omega),$$
  and then for all $k \in \mathbb{N} \cup \{0\}$ we have a canonical $\mathcal{O}$-module isomorphism
  $$\mathscr{S}^k_{\mathcal{O}}\mathcal{L} \cong \mathscr{U}_k(\mathcal{O}, \mathcal{L}_{\omega}, t)/{\mathscr{U}_{k-1}(\mathcal{O}, \mathcal{L}_{\omega}, t)}.$$
   At the level of associated presheaves of $\mathscr{U}(\mathcal{O}, \mathcal{L}, \omega)$ and $\mathscr{S}_{\mathcal{O}}\mathcal{L}$, etc., the description follows from the Lie–Rinehart algebra case (see \cite[Theorem 3.3 and Lemma 3.6]{HM}). More precisely, 
for a point $x \in X$, there exists an open set $U_x \subset X$, such that $\mathcal{L}|_{U_x}$ is a free $\mathcal{O}|_{U_x}$-module. Thus, for any open subset $U \subset U_x$, for the free Lie-Rinehart algebra $(\mathcal{O}(U), \mathcal{L}(U))$ together with the $2$-cocycle $w_U$, we have canonical (local) isomorphisms 
   $$\mathcal{U}(\mathcal{O}(U), \mathcal{L}(U)_{w_U}, t_U) \cong \mathcal{U}(\mathcal{O}(U), \mathcal{L}(U), w_U),~~{S}^k_{\mathcal{O}(U)} \mathcal{L}(U) \cong \mathcal{U}_k(\mathcal{O}(U), \mathcal{L}(U)_{w_U}, t_U)/{\mathcal{U}_{k-1}(\mathcal{O}(U), \mathcal{L}(U)_{w_U}, t_U)},$$ of $\mathbb{K}$-algebras and of  $\mathcal{O}(U)$-modules,
 respectively. Both these isomorphisms are compatible with the associated restriction maps.
   Extending these constructions from the local to its global setting by applying sheafification and hypercohomology functor, thereby obtaining the desired result.
\end{proof}
\begin{Cor}
    In particular, for a Lie algebroid $\mathcal{L}$ over a complex manifold $X$ which is a locally free $\mathcal{O}_X$-module of finite rank $($i.e. $\mathcal{L}$ is a holomorphic Lie algebroid$)$, we get the result given in \cite{PT}.
\end{Cor}
\subsubsection{\textbf{Applications to particular cases}}	In the following we describe Theorem \ref{connection with prescribed curvature} and Theorem \ref{twisted PBW theorem} in the Lie-Rinehart algebras and (smooth or holomorphic) Lie algebroids context.

\medskip

\noindent $(1)$ \textbf{On Lie-Rinehart algebras.} Let $(R, L)$ be a Lie-Rinehart algebra, with projective left $R$-module structure on $L$. Also, let $w$ be a representative of a cohomology class in the second Lie-Rinehart cohomology group $H^2(L, R)$. This is equivalent to considering a locally free (quasicoherent) Lie algebroid $\mathcal{L}$ over $(X, \mathcal{O}_X)$ with a cohomology class $\bar{\omega} \in \mathbb{H}^2(\mathcal{L}, \mathcal{O}_X)$ where $X=Spec(R)$.  Notice that, the twisted universal enveloping algebra $\mathcal{U}(R, L, w)$ is canonically isomorphic to the space of global sections $\mathscr{U}(\mathcal{O}_X, \mathcal{L}, \omega)(X)$. Thus, using Theorem \ref{connection with prescribed curvature}, an $R$-module $E$ with an $L$-connection of curvature type $w$ (see \cite{JH}) is equivalent to a left $\mathcal{U}(R, L, w)$-module structure on $E$. In this situation, Theorem \ref{twisted PBW theorem} provides a canonical isomorphism  of graded $R$-algebras between the graded quotient $gr(\mathcal{U}(R, L, w))$ of the almost commutative algebra $\mathcal{U}(R, L, w)$ and the symmetric algebra $S_RL$ of the $R$-module $L$. 

\medskip

\noindent $(2)$ \textbf{On holomorphic Lie algebroids.} For a holomorphic Lie algebroid \( L \) over a complex manifold \( X \), the existence of a flat \( L \)-connection, or even a general \( L \)-connection on a holomorphic vector bundle $E$ is not guaranteed. Equivalently, consider the locally free \( \mathcal{O}_X \)-module \(\Gamma_X(E)=: \mathcal{E} \)   of finite rank, and studying 
 $\Gamma_X(L)=:\mathcal{L}$-connections on \( \mathcal{E} \).
In particular, when $\mathcal{L}=\mathcal{T}_X$, a $\mathcal{T}_X$-connection on $\mathcal{E}$ (i.e., a usual holomorphic connection on $E$) is
a global section \( \nabla : \mathcal{T}_X \to \mathcal{A}t(\mathcal{E}) \) of the anchor $\sigma: \mathcal{A}t(\mathcal{E}) \rightarrow \mathcal{T}_X$ satisfying \( \sigma \circ \nabla = \mathrm{Id}_{\mathcal{T}_X} \) (see ~\eqref{L-connection}), provides a splitting of the short exact sequence~\eqref{Atiyah sequence} of \( \mathcal{O}_X \)-modules. If the curvature of the connection \( \nabla \) vanishes, the sequence splits as Lie algebroids.

In contrast, in the smooth category, Lie algebroid connections on smooth vector bundles always exist due to the availability of partitions of unity. Over certain complex manifolds, such as Stein manifolds, a holomorphic Lie algebroid connection always exists (see~\cite{Atiyah}, \cite{PT}). On a Riemann surface \( X \), a \( \mathcal{T}_X \)-connection on \( \mathcal{E} \) corresponds to a flat holomorphic connection in the classical sense, exists (see~\cite{AS-IB}).

These considerations motivate the study of \( \mathcal{L} \)-connections with nonzero curvature on a locally free sheaf \( \mathcal{E} \), where \( \mathcal{L} \) is a locally free Lie algebroid over \( (X, \mathcal{O}_X) \). By Theorem~\ref{connection with prescribed curvature}, such connections with prescribed curvature class \( \bar{\omega} \in \mathbb{H}^2(\mathcal{L}, \mathcal{O}_X) \) correspond to modules over the sheaf of filtered algebras $\mathscr{U}(\mathcal{O}_X, \mathcal{L}, \omega)$. This provides the result 
[Proposition $5.3$, \cite{PT}], with relaxing the coherent or finite rank condition on $\mathcal{L}$ over $\mathcal{O}_X$. 

\medskip

\noindent $(3)$ \textbf{Logarithmic foliations and Meromorphic connections.} 
Let \((Y, \mathcal{O}_Y)\) be an analytic subspace (or a divisor) in a complex manifold \((X, \mathcal{O}_X)\). Consider the Lie algebroid given by the sheaf of logarithmic vector fields along \(Y\), i.e. \(\mathcal{L} := \mathcal{T}_X(-\log~ Y)\). An \(\mathcal{L}\)-connection on an \(\mathcal{O}_X\)-module \(\mathcal{E}\) corresponds to a meromorphic connection on \(\mathcal{E}\) with simple poles along \(Y\) (see \cite{BP}, \cite{Simpson}) is described by a $\mathbb{C}_X$-linear map
$$\nabla : \mathcal{E} \rightarrow \Omega^1_{X}(log~Y) \otimes_{\mathcal{O}_X}\mathcal{E},$$
satisfying the Leibniz rule (see (\ref{L-connection special})). In \cite{Gualtieri-Pym}, such Lie algebroids are used to study meromorphic connections and their relationship with differential equations, leading to a novel method for studying divergent asymptotic series near irregular singularities.
Such \(\mathcal{L}\)-connections on \(\mathcal{E}\) can also be studied via the theory of modules over a twisted version of the sheaf of logarithmic differential operators \(\mathcal{D}_X(-\log~ Y)\).

\medskip

\noindent $(4)$ \textbf{Poisson vector bundles and contravariant connections.} In the smooth category, Ginzburg introduced the notion of Poisson vector bundles~\cite{VGinz1}, which was later extended by Bursztyn to Poisson modules~\cite{HB}, thereby generalizing the classical notions of Poisson manifolds and Poisson algebras, respectively. In particular, it is shown that a contravariant connection on a vector bundle corresponds to a Poisson vector bundle structure on it (see [Proposition~2.8, ~\cite{HB}]).
This relationship can be reformulated as follows: \\
Let $(X, \mathcal{O}_X)$ be a (smooth or holomorphic) Poisson manifold, and let $\mathcal{E}$ be a locally free $\mathcal{O}_X$-module of finite rank. The structure sheaf $\mathcal{O}_X$ carries a canonical Poisson algebra structure, and $\Omega^1_X$ inherits a natural Lie algebroid structure (see Example~\ref{Cotantgent Poisson}). A contravariant connection on $\mathcal{E}$, which is dual to the usual \( \mathcal{T}_X \)-connection, i.e. dual of the notion as described in~\eqref{L-connection special}, as an action (see \cite{RF}) given by
\[
\nabla: \Omega^1_X \times \mathcal{E} \rightarrow \mathcal{E}
\]
satisfying standard properties: $\mathcal{O}_X$-linearity in the first argument and the Leibniz rule for the $\mathcal{O}_X$-action on $\mathcal{E}$.
A Poisson module structure on $\mathcal{E}$ is then equivalent to the existence of a flat contravariant connection on $\mathcal{E}$, i.e., a (Lie algebroid) \( \Omega^1_X \)-module structure. Therefore, to study Poisson module structures on (smooth or holomorphic) vector bundles, it suffices to work with the notion of \( \mathscr{U}(\mathcal{O}_X, \Omega^1_X) \)-module structures. More generally, using the notion of Poisson cohomology (see \cite{HB, LG-PX, JH}), one obtains a twisted analogue.

\subsection{Deformation groupoid of a Lie algebroid.}
   We introduce the notion of the deformation groupoid of a Lie algebroid over a ringed space. This construction provides a global counterpart to the deformation groupoid theory of Lie-Rinehart algebras as studied in \cite{HM}. This framework is also related to the moduli space of sheaves of rings of differential operators satisfying Simpson's axioms \cite{PT}. 
   
   Let  $\mathcal{A}(\mathscr{S}_{\mathcal{O}} \mathcal{L})$ be the category, whose objects are of the form $(\mathscr{U}, \psi_{\mathscr{U}})$, where $\mathscr{U}$ is a almost commutative filtered algebra over $(X, \mathcal{O})$ and is equipped with an isomorphism of graded $\mathcal{O}$-algebras
   \begin{equation} \label{almost commutative}
      \psi_{\mathscr{U}}: \mathscr{S}_{\mathcal{O}} \mathcal{L} \rightarrow gr(\mathscr{U}).
   \end{equation}
    A morphism 
   $$\phi: (\mathscr{U}, \psi_{\mathscr{U}}) \rightarrow (\mathscr{U}', \psi_{\mathscr{U}'})$$ in $\mathcal{A}(\mathcal{S}_{\mathcal{O}} \mathcal{L})$ is a sheaf homomorphism of filtered $\mathbb{K}_X$-algebras $\phi: \mathscr{U} \rightarrow \mathscr{U}'$ satisfying $gr(\phi)  \circ \psi_{\mathscr{U}} = \psi_{\mathscr{U}'}$, where $gr(\phi)$ is the induced map on the  associated sheaves of $\mathbb{K}_X$-algebras. Since $\psi_{\mathscr{U}}$ and $\psi_{\mathscr{U}'}$ are isomorphisms, the associated graded map 
   $$gr(\phi)= \psi_{\mathscr{U}'} \circ \psi_{\mathscr{U}}^{-1}$$ is an isomorphism of graded $\mathcal{O}$-algebras. It follows that, the map $\phi$ is an isomorphism of filtered $\mathbb{K}_X$-algebras over $(X, \mathcal{O})$. In this category, morphisms are invertible since they are isomorphisms, and thus this category forms a groupoid. Also, this groupoid captures information for the deformation of the classical graded filtered algebra $\mathscr{S}_{\mathcal{O}}\mathcal{L}$. We refer the category $\mathcal{A}(\mathscr{S}_{\mathcal{O}} \mathcal{L})$ by deformation groupoid of a Lie algebroid $(\mathcal{O},\mathcal{L})$ over $X$.

\begin{Rem}
When \( (X, \mathcal{O}_X) \) is a smooth algebraic variety (or a smooth manifold, or a complex manifold), the following identifications hold. If \( \mathscr{U} = \mathcal{D}_X \), then the category of \( \mathscr{U} \)-modules corresponds to the category of (usual) flat connections. If \( \mathscr{U} = \mathscr{S}_{\mathcal{O}_X} \mathcal{T}_X \), then semistable \( \mathscr{U} \)-modules are equivalent to semistable Higgs sheaves. More generally, \( \mathscr{U} \) can be viewed as a one-parameter deformation of \( \mathscr{S}_{\mathcal{O}_X} \mathcal{T}_X \) to \( \mathcal{D}_X \). For further details, we refer the reader to Simpson's foundational work \cite{Simpson}. A similar framework can be extended to the setting of Lie algebroids, providing a natural context to understand deformation groupoids; see \cite{PT, HM}.
\end{Rem}

 We extend some of the important results described in \cite{HM} for the deformation groupoid as follows. First, we relate the cohomology group $\mathbb{H}^2(\mathcal{L}, \mathcal{O})$ with deformations groupoid $\mathcal{A}(\mathscr{S}_{\mathcal{O}} \mathcal{L})$ of the sheaf of almost commutative filtered algebras  $(\mathscr{U}, \mathscr{U}_i)$ where $\mathcal{L} \cong \mathscr{U}_1/{\mathscr{U}_0}$ is a locally free Lie algebroid.  We show that there is a one-to-one correspondence between
	$\mathbb{H}^2(\mathcal{L},\mathcal{O})$ and the set of isomorphic classes of objects in $\mathcal{A}(\mathscr{S}_{\mathcal{O}} \mathcal{L})$.

\begin{Rem} \label{U & w}
    The isomorphism (\ref{almost commutative}) for an object $(\mathscr{U}, \mathscr{U}_i)$ in the deformation groupoid $\mathcal{A}(\mathscr{S}_{\mathcal{O}} \mathcal{L})$ 
 describes an isomorphism $\mathscr{U}_1/{\mathcal{O}} \cong \mathcal{L}$, which also splits as Lie algebroids $\mathscr{U}_1 \cong \mathcal{O} \oplus \mathcal{L}$ when $\mathcal{L}$ is locally free $\mathcal{O}$-module. This provides an abelian extension of $\mathcal{L}$ by $\mathcal{O}$, and thus (\ref{Ext-H^2}) corresponds to a cohomology class $\bar{\omega} \in \mathbb{H}^2(\mathcal{L},\mathcal{O})$.
\end{Rem}
Let  $Iso(\mathcal{A}(\mathscr{S}_{\mathcal{O}} \mathcal{L}))$ be the set of isomorphism classes of objects in $\mathcal{A}(\mathscr{S}_{\mathcal{O}} \mathcal{L})$. Define the following map
$$\psi: \mathbb{H}^2(\mathcal{L},\mathcal{O}) \rightarrow Iso(\mathcal{A}(\mathscr{S}_{\mathcal{O}} \mathcal{L}))$$
given by 
\begin{equation} \label{cohomology-deformation grouoid}
    \psi(\bar{\omega}):= (\mathscr{U}(\mathcal{O}, \mathcal{L}, \omega), \phi_{\omega}),
\end{equation}
where $\phi_{\omega}: \mathscr{S}_{\mathcal{O}}\mathcal{L} \rightarrow gr(\mathscr{U}(\mathcal{O}, \mathcal{L}, \omega)) $ is the canonical isomorphism. The map $\psi$ is well defined since for two cocycles $\omega$ and $\omega+ d^1 \rho$, which are differs by a coboundary, representing the same cohomology class $\bar{\omega} \in \mathbb{H}^2(\mathcal{L},\mathcal{O})$ we can show there is an isomorphism
$$\mathscr{U}(\mathcal{O}, \mathcal{L}, \omega) \cong \mathscr{U}(\mathcal{O}, \mathcal{L}, \omega+ d^1 \rho)$$
of sheaf of filtered algebras in $\mathcal{A}(\mathcal{S}_{\mathcal{O}} \mathcal{L})$.

Now we state and prove the above result in an explicit form as follows.
\begin{Thm} \label{deformation groupoid}
Let $\mathcal{L}$ be a locally free Lie algebroid over $(X, \mathcal{O})$. The map $\psi: \mathbb{H}^2(\mathcal{L},\mathcal{O}) \rightarrow Iso(\mathcal{A}(\mathscr{S}_{\mathcal{O}} \mathcal{L}))$ defined in (\ref{cohomology-deformation grouoid})
provides a one-to-one correspondence of sets.
\end{Thm}
\begin{proof} To show the one-to-one property of the map $\psi$, observe that an isomorphism
$$\mathscr{U}(\mathcal{O}, \mathcal{L}, \omega) \cong \mathscr{U}(\mathcal{O}, \mathcal{L}, \omega')$$ implies the associated cohomology classes $\bar{\omega}=\bar{\omega}' \in  \mathbb{H}^2(\mathcal{L},\mathcal{O})$ (i.e. $\omega$ and $\omega'$ are cohomologous).

We now show that the map $\psi$  is surjective. Let $(\mathscr{U}, \mathscr{U}_i)$ be an isomorphism class of objects in $\mathcal{A}(\mathcal{S}_{\mathcal{O}} \mathcal{L})$. Thus, we need to show there exists a twisted universal enveloping algebroid $\mathscr{U}(\mathcal{O}, \mathcal{L}, \omega)$ isomorphic to $\mathscr{U}$ for some $2$-cocycle $\omega$ associated to $ \mathbb{H}^2(\mathcal{L},\mathcal{O})$.
In Remark \ref{U & w}, we mention that how $\mathscr{U}_1$ is associated to the Lie algebroid $\mathcal{L}$ and a cohomology class $\bar{w} \in \mathbb{H}^2(\mathcal{L}, \mathcal{O})$. Now, we show that $\mathscr{U} \cong \mathscr{U}(\mathcal{O}, \mathcal{L}, \omega)$ as filtered associative algebras. 
From the construction described in Section \ref{Twisted enveloping algebroid}, we have a canonical map 
  $$\bar{\mathcal{P}}: \mathscr{T}^1(\mathcal{L}_{\omega}) \rightarrow  \mathscr{U}(\mathcal{O}, \mathcal{L}, \omega).$$
  Using the locally free condition on $\mathcal{L}$, consider the map $\tilde{\psi}: \mathscr{T}^1(\mathcal{L}_{\omega}) \rightarrow  \mathscr{U}$ defined by
  $$\tilde{\psi}((f_1 s+D_1) \otimes \cdots \otimes (f_k s+D_k)):= \textstyle\prod_i(f_i+ (\alpha \circ \psi_{\mathcal{L}})(D_i)),$$
  where $\psi_{\mathcal{L}}:=\psi_{\mathscr{U}}|_{\mathcal{L}}$ is the restriction of the graded algebra isomorphism (\ref{almost commutative}) on $\mathcal{L}$ and there is a section $\alpha: \mathscr{U}_1/{\mathcal{O}} \rightarrow \mathscr{U}_1$ of the projection map $p: \mathscr{U}_1 \rightarrow  \mathscr{U}_1/{\mathcal{O}}$ associated with the abelian Lie algebroid extension
  $$0 \rightarrow \mathcal{O} \rightarrow \mathscr{U}_1 \overset{p}{\rightarrow} \mathscr{U}_1/{\mathcal{O}} \rightarrow 0.$$
  Notice that, the map $\tilde{\psi}$ is an $\mathcal{O}$-module homomorphism and satisfies
  $$\tilde{\psi}((f_1 s+D_1) \otimes (f_2 s+D_2)- (f_2 s+D_2) \otimes (f_1 s+D_1)- [f_1 s+D_1, f_2 s+D_2]_{\omega})=0.$$ 
  Thus, the map $\tilde{\psi}$ induces a sheaf of filtered algebra isomorphism 
  $$\psi': \mathscr{U}(\mathcal{O}, \mathcal{L}, \omega) \rightarrow \mathscr{U}$$
  defined by 
  $$\psi'(D_1 \cdots D_k)=(\alpha \circ \psi_{\mathcal{L}})(D_1) \cdots (\alpha \circ \psi_{\mathcal{L}})(D_k).$$
  For the underlying presheaves level (or local) descriptions, observe that for each open set $U$ of $X$, we have $\mathcal{O}(U)$-linear $\mathbb{K}$-algebra homomorphism
  $\psi'_U: \mathcal{U}(\mathcal{O}(U), \mathcal{L}(U), w) \rightarrow \mathcal{U}(U)$. The homomorphisms satisfies the compatibility condition $res_{UV}^{\mathcal{U}} \circ \psi'_U= \psi'_V \circ res_{UV}^{\mathcal{U}(\mathcal{O}, \mathcal{L}, w)}$, for any open subset $V$ of $U$. Since the Lie algebroid $\mathcal{L}$ is locally free $\mathcal{O}$-module, for a point $x \in X$, there exists an open set $U_x \subset X$ around $x$, such that $\mathcal{L}|_{U_x}$ is a free $\mathcal{O}|_{U_x}$-module. Thus, for any open subset $U \subset U_x$, we have canonical  isomorphisms $\psi'_U$.
   Hence,  using the sheafication of the presheaf homomorphism $U \mapsto \psi'_U$, we have the above isomorphism $\psi'$. For details description in the projective Lie–Rinehart algebra case, see \cite[Lemmas 4.11 and 4.12]{HM}.
\end{proof}

\begin{Cor}
 Assume $\mathcal{L}:= \mathscr{U}_1/{\mathscr{U}_0}$ is a locally free $\mathcal{O}$-module. Therefore, the the Lie algebroid $\mathcal{L}$ together with a $2$-cocycle $\omega \in \mathcal{Z}^2(\mathcal{L}, \mathcal{O})$ provides an isomorphism $\mathscr{U} \cong \mathscr{U}(\mathcal{O}, \mathcal{L}, \omega)$ of filtered associative algebras.
\end{Cor}
\begin{Rem}
    In the category $\mathcal{A}(\mathscr{S}_{\mathcal{O}} \mathcal{L})$, morphisms are parametrized by the sheaf of $1$-cocycles $\mathcal{Z}^1(\mathcal{L}, \mathcal{O})$ of the cochain complex (\ref{Chevalley-Eilenberg-de Rham complex}) for the Lie algebroid $(\mathcal{O}, \mathcal{L})$ with coefficient in $\mathcal{O}$ (see \cite{HM} for the local descriptions).
    \end{Rem}

\subsubsection{\textbf{Application to particular cases}.} We describe existing notions related to deformation groupoids or moduli spaces for Lie algebroid connections as follows.

\medskip

\noindent $(1)$ \textbf{Lie-Rinehart algebras setting.} Recall that a locally free Lie algebroid $\mathcal{L}$ (of finite rank) over an affine scheme $(X, \mathcal{O}_X)$ corresponds to the (finitely generated) projective Lie-Rinehart algebra $L := \mathcal{L}(X)$ over the ring $R := \mathcal{O}_X(X)$, and vice versa. Thus, we recover the algebraic description of the deformation groupoid \( A(\mathrm{Sym}^*_R(L)) \) of a Lie-Rinehart algebra \( L \), projective as a left \( R \)-module, as studied by Maakestad in \cite{HM}. 
The isomorphism classes in this groupoid are parametrized by \( H^2(L, R) \), and morphisms by \( Z^1(L, R) \). 
As an application, for any filtered almost commutative ring \( \mathcal{U} \) satisfying the PBW condition  $gr(\mathcal{U}) \cong S_R L$, we have the category of left \( \mathcal{U} \)-modules is equivalent to the category of \( L \)-connections with curvature type \( w \in Z^2(L, R) \), where $\mathcal{U} \cong \mathcal{U}(R, L, w)$. This generalizes the theory of deformation groupoid originally introduced by Sridharan for Lie algebras that are free modules over a commutative algebra \cite{Sridh}.

In particular, for $R= C^{\infty}(M)$ and $L$ is the space of global sections of a smooth Lie algebroid over a smooth manifold $M$, we get analogous results in the context of differential geometry. In particular, if $L= Der_{\mathbb{R}}(C^{\infty}(M))$, then we recover the algebras of twisted differential operators over $M$.

\medskip

\noindent $(2)$ \textbf{Holomorphic Lie algebroids setting.} 
   Let \( (\mathscr{U}, \psi_{\mathscr{U}}) \) be a sheaf of almost commutative filtered \( \mathbb{C} \)-algebras over a complex manifold \( X \) satisfying Simpson’s axioms~\eqref{Simpson's axioms}.
   Assume further that \( \mathcal{L} := \mathscr{U}_1 / \mathscr{U}_0 \) is a locally free \( \mathcal{O}_X \)-module. Thus,  \( (\mathscr{U}, \psi_{\mathscr{U}}) \in  \mathcal{A}(\mathscr{S}_{\mathcal{O}} \mathcal{L})\) and \( \mathcal{L} \) naturally acquires the structure of a holomorphic Lie algebroid, and \( \mathscr{U}_1 \) becomes an abelian Lie algebroid extension of \( \mathcal{L} \) by \( \mathcal{O}_X \cong \mathscr{U}_0 \), classified by a class \( \bar{\omega} \in \mathbb{H}^2(\mathcal{L}, \mathcal{O}_X) \). Conversely, for a holomorphic Lie algebroid \( \mathcal{L} \) over \( (X, \mathcal{O}_X) \) and a class \( \bar{\omega} \in \mathbb{H}^2(\mathcal{L}, \mathcal{O}_X) \), the associated twisted universal enveloping algebroid \( \mathscr{U} := \mathscr{U}(\mathcal{O}_X, \mathcal{L}, \omega) \) satisfies Simpson’s axioms~\eqref{Simpson's axioms}.

Over a compact Kähler manifold (or smooth complex projective variety) $(X, \mathcal{O}_X)$, a stronger formulation of this correspondence has been established by Tortella in~\cite{PT},  as follows.
 There is a canonical one-to-one correspondence between pairs \( (\mathscr{U}, \psi_{\mathscr{U}}) \), with \( \mathscr{U} \) a sheaf of almost polynomial filtered algebras over $(X, \mathcal{O}_X)$ satisfying Simpson’s axioms, and pairs \( (\mathcal{L}, \omega) \), where \( \mathcal{L} \) is a holomorphic Lie algebroid structure on \( \mathscr{U}_1 / \mathscr{U}_0 \), and \( \omega \in F^1 H^2(\mathcal{L}, \mathbb{C}) \) lies in the first Hodge filtration of the Lie algebroid cohomology.

\medskip

\noindent $(3)$ \textbf{Logarithmic foliation and (logarithmic) de Rham cohomology.}  For $\mathcal{L}= \mathcal{T}_X$ of a non-singular ringed space, we have the (smooth, holomorphic and algebraic) de Rham theorem, due to de Rham, Serre and Grothendieck, respectively. When $X$ is a complex manifold or a smooth scheme over $\mathbb{C}$, it establishes an isomorphism between the de Rham hypercohomology \( \mathbb{H}^{\bullet}(X, \Omega^{\bullet}_X) \), the sheaf cohomology $H^{\bullet}(X, \mathbb{C}_X)$ and the singular cohomology \( H_{\mathrm{sing}}^{\bullet}(X(\mathbb{C}), \mathbb{C}) \), where $X(\mathbb{C})$ is the analytification of $X$ (see \cite{VID}).
Moreover, for  a divisor $Y$ of $X$, the Chevalley-Eilenberg-de Rham complex of $\mathcal{L}=\mathcal{T}_X(-log~Y)$ is  the logarithmic de Rham complex of meromorphic forms $\Omega^{\bullet}_X(log ~Y)$ over $X$ with simple poles along the divisor $Y$. The associated  hypercohomology  $\mathbb{H}^{\bullet}(X, \Omega^{\bullet}_X(log ~Y))$ is isomorphic to the singular cohomology $H_{sing}^{\bullet}(X \setminus Y, \mathbb{C})$ for a normal crossing divisor $Y$ of a smooth complex projective variety $X$ (see \cite{PT, AS}). 
Using these cohomological results in smooth, holomorphic, and algebraic settings, one can simplify the study of the deformation groupoid for the Lie algebroids of the tangent sheaf $\mathcal{T}_X$ and the logarithmic tangent sheaf $\mathcal{T}_X(-log~Y)$.

    \begin{Rem}
In Example \ref{free Lie algebroid}, the universal enveloping algebroid $\mathbb{D}_X$ of the path algebroid $\mathcal{P}_X$ of a smooth ringed space $(X, \mathcal{O}_X)$ is recalled. A $\mathcal{P}_X$-module is simply a vector bundle with a (not necessarily flat) $\mathcal{T}_X$-connection, studied as a $\mathbb{D}_X$-module, describing systems of linear PDEs on the path space (see \cite{MK}).
\end{Rem}
\begin{Rem}
    $($\emph{Universal algebra.}$)$ In [Theorem A.13, \cite{HM}], Maakestad introduced the notion of universal algebra $U^{ua}(L)$ of a Lie–Rinehart algebra $L$ $($over an algebra$)$ and showed the category of $U^{ua}(L)$-modules is equivalent to the category of $L$-connections. This categorical equivalence can be extended to the algebro-geometric setting in order to study the category of $\mathcal{L}$-connections for a Lie algebroid $\mathcal{L}$ over a ringed space.
\end{Rem}
\subsection{Remarks.} 

\emph{$\check{C}$ech Cohomology.} For a non-singular ringed space $(X, \mathcal{O}_X)$, one can find good (e.g., affine or Stein) open covers of $X$  that allow standard computations of sheaf cohomology using Čech cohomology. This method naturally extends to the computation of hypercohomology for complexes of $\mathcal{O}_X$-modules (see \cite{PT, BRT, BP, AS}). A classical example of this approach is the Čech–de Rham cohomology. Therefore, we may use Čech hypercohomology to simplify the computations and formulations of the above results.

\emph{Stalkwise description.} All our sheaf-theoretic constructions can be understood in terms of the induced algebraic structures on the stalks. Since the stalks of a presheaf and its associated sheaf coincide, the descriptions involving presheaves or local data are equivalent to those given at the level of stalks (see \cite{SR, UB}).

\section*{Conclusion}
In this article, we studied Lie algebroids, quantum Poisson algebroids, and Lie algebroid connections. The results in this paper also leads to the following directions:  

$(1)$ \textbf{Automorphisms.} The study of automorphisms of the quantum Poisson algebra $\mathcal{D}(M)$ offers deep geometric insight of a smooth manifold $M$ \cite{GP}. We aim to investigate how automorphisms of the universal enveloping algebroid $\mathscr{U}(\mathcal{O}, \mathcal{L})$ reflect the geometry of a Lie algebroid $(\mathcal{O}, \mathcal{L})$, extending classical results from \cite{GP} to Lie algebroids such as tangent sheaf, sheaf of logarithmic derivations and Atiyah algebroids, etc.

$(2)$ \textbf{Cohomology of quantum Poisson algebras.} It would be interesting to develop a cohomology theory for quantum Poisson algebras. Such a theory could capture deformation-theoretic and representation-theoretic aspects intrinsic to filtered algebras and noncommutative geometries.

$(3)$ \textbf{Over ringed sites.} The notion of a ringed site generalizes that of a ringed space. Certain aspects of the theory of smooth Lie algebroids have been extended to the broader setting of Lie algebroids over ringed sites (see \cite{DRV, CRV, CV}). Accordingly, parts of the theory developed in this article may also extend to that context.

\vspace{.2 cm}
\textbf{Acknowledgments.} 
 The authors express their heartfelt gratitude to Dr. Ashis Mandal for introducing them to several foundational concepts that have significantly shaped this work. The second author is grateful to Prof. Mainak Poddar for some of his valuable suggestions and acknowledges the support of IISER Pune for the Institute Post-Doctoral fellowship IISER-P/Jng./20235445.
 
\vspace{.1 cm}

\bibliographystyle{plain}

\begin{thebibliography}{10}

\bibitem{Atiyah}
M.~F. Atiyah.
\newblock Complex analytic connections in fibre bundles.
\newblock {\em Trans. Amer. Math. Soc.}, 85:181--207, 1957.

\bibitem{NK-PS}
Xavier Bekaert, Niels Kowalzig, and Paolo Saracco.
\newblock Universal enveloping algebras of {L}ie--{R}inehart algebras: crossed products, connections, and curvature.
\newblock {\em Lett. Math. Phys.}, 114(6):Paper No. 140, 2024.

\bibitem{AS-IB}
Indranil Biswas and Anoop Singh.
\newblock Line bundles on the moduli space of {L}ie algebroid connections over a curve.
\newblock {\em Bull. Sci. Math.}, 193:Paper No. 103421, 30, 2024.

\bibitem{UB}
Ugo Bruzzo.
\newblock Lie algebroid cohomology as a derived functor.
\newblock {\em J. Algebra}, 483:245--261, 2017.

\bibitem{BRT}
Ugo Bruzzo, Igor Mencattini, Vladimir~N. Rubtsov, and Pietro Tortella.
\newblock Nonabelian holomorphic {L}ie algebroid extensions.
\newblock {\em Internat. J. Math.}, 26(5):1550040, 26, 2015.

\bibitem{HB}
Henrique Bursztyn.
\newblock Poisson vector bundles, contravariant connections and deformations.
\newblock {\em Progr. Theoret. Phys. Suppl.}, Number 144, pages 26--37. 2001.


\bibitem{DRV}
Damien Calaque, Carlo~A. Rossi, and Michel van~den Bergh.
\newblock Hochschild (co)homology for {L}ie algebroids.
\newblock {\em Int. Math. Res. Not. IMRN}, (21):4098--4136, 2010.

\bibitem{CRV}
Damien Calaque, Carlo~A. Rossi, and Michel Van~den Bergh.
\newblock C\u ald\u araru's conjecture and {T}sygan's formality.
\newblock {\em Ann. of Math. (2)}, 176(2):865--923, 2012.

\bibitem{CV}
Damien Calaque and Michel Van~den Bergh.
\newblock Hochschild cohomology and {A}tiyah classes.
\newblock {\em Adv. Math.}, 224(5):1839--1889, 2010.

\bibitem{Francisco}
Francisco~J. Calder\'on-Moreno.
\newblock Logarithmic differential operators and logarithmic de {R}ham complexes relative to a free divisor.
\newblock {\em Ann. Sci. \'Ecole Norm. Sup. (4)}, 32(5):701--714, 1999.

\bibitem{CMD}
Francisco~J. Castro-Jim\'enez, Luis Narv\'aez-Macarro, and David Mond.
\newblock Cohomology of the complement of a free divisor.
\newblock {\em Trans. Amer. Math. Soc.}, 348(8):3037--3049, 1996.

\bibitem{SC}
Sophie Chemla.
\newblock A duality property for complex {L}ie algebroids.
\newblock {\em Math. Z.}, 232(2):367--388, 1999.

\bibitem{VID}
V.~I. Danilov.
\newblock Cohomology of algebraic varieties.
\newblock In {\em Current problems in mathematics. {F}undamental directions, {V}ol. 35 ({R}ussian)}, Itogi Nauki i Tekhniki, pages 5--130, 272. Akad. Nauk SSSR, Vsesoyuz. Inst. Nauchn. i Tekhn. Inform., Moscow, 1989.

\bibitem{Jong-Johan-Max-Shin}
Aise~Johan de~Jong, Max Lieblich, and Minseon Shin.
\newblock Locally free twisted sheaves of infinite rank.
\newblock {\em Doc. Math.}, 28(1):133--171, 2023.

\bibitem{RF}
Rui~Loja Fernandes.
\newblock Lie algebroids, holonomy and characteristic classes.
\newblock {\em Adv. Math.}, 170(1):119--179, 2002.

\bibitem{VGinz2}
Viktor~L. Ginzburg.
\newblock Lectures on $\mathcal{D}$-modules.
\newblock {\em lecture notes of a course given at the University of Chicago}, 1998.

\bibitem{VGinz1}
Viktor~L. Ginzburg.
\newblock Grothendieck groups of {P}oisson vector bundles.
\newblock {\em J. Symplectic Geom.}, 1(1):121--169, 2001.

\bibitem{GuTw}
Ulrich G\"ortz and Torsten Wedhorn.
\newblock {\em Algebraic geometry {II}: {C}ohomology of schemes---with examples and exercises}.
\newblock Springer Studium Mathematik---Master. Springer Spektrum, Wiesbaden, [2023] \copyright 2023.

\bibitem{GP}
J.~Grabowski and N.~Poncin.
\newblock Automorphisms of quantum and classical {P}oisson algebras.
\newblock {\em Compos. Math.}, 140(2):511--527, 2004.

\bibitem{Gualtieri-Pym}
Marco Gualtieri, Songhao Li, and Brent Pym.
\newblock The {S}tokes groupoids.
\newblock {\em J. Reine Angew. Math.}, 739:81--119, 2018.

\bibitem{JH}
Johannes Huebschmann.
\newblock Poisson cohomology and quantization.
\newblock {\em J. Reine Angew. Math.}, 408:57--113, 1990.

\bibitem{MK}
Mikhail Kapranov.
\newblock Free {L}ie algebroids and the space of paths.
\newblock {\em Selecta Math. (N.S.)}, 13(2):277--319, 2007.

\bibitem{KM-KS}
Y.~Kosmann-Schwarzbach and K.~C.~H. Mackenzie.
\newblock Differential operators and actions of {L}ie algebroids.
\newblock In {\em Quantization, {P}oisson brackets and beyond ({M}anchester, 2001)}, volume 315 of {\em Contemp. Math.}, pages 213--233. Amer. Math. Soc., Providence, RI, 2002.

\bibitem{LG}
Camille Laurent-Gengoux and Leonid Ryvkin.
\newblock The neighborhood of a singular leaf.
\newblock {\em J. \'Ec. polytech. Math.}, 8:1037--1064, 2021.

\bibitem{LG-PX}
Camille Laurent-Gengoux, Mathieu Sti\'enon, and Ping Xu.
\newblock Holomorphic {P}oisson manifolds and holomorphic {L}ie algebroids.
\newblock {\em Int. Math. Res. Not. IMRN}, pages Art. ID rnn 088, 46, 2008.

\bibitem{HM-Chernclass}
Helge Maakestad.
\newblock The {C}hern character for {L}ie-{R}inehart algebras.
\newblock {\em Ann. Inst. Fourier (Grenoble)}, 55(7):2551--2574, 2005.

\bibitem{HM}
Helge~{\O}ystein Maakestad.
\newblock Algebraic connections on projective modules with prescribed curvature.
\newblock {\em J. Algebra}, 436:161--227, 2015.

\bibitem{LM}
Luis~Narv$\acute{a}$ez Macarro.
\newblock Differential structures in commutative algebra.
\newblock {\em Mini-course at the XXIII Brazilian Algebra Meeting, July 27 - August 1}, Maring$\acute{a}$, Brazil., 2014.

\bibitem{KM}
Kirill C.~H. Mackenzie.
\newblock {\em General theory of {L}ie groupoids and {L}ie algebroids}, volume 213 of {\em London Mathematical Society Lecture Note Series}.
\newblock Cambridge University Press, Cambridge, 2005.

\bibitem{AA}
Ashis Mandal and Abhishek Sarkar.
\newblock On {L}ie algebroid over algebraic spaces.
\newblock {\em Comm. Algebra}, 51(4):1594--1613, 2023.

\bibitem{MM}
I.~Moerdijk and J.~Mr\v{c}un.
\newblock On the universal enveloping algebra of a {L}ie algebroid.
\newblock {\em Proc. Amer. Math. Soc.}, 138(9):3135--3145, 2010.

\bibitem{JNes}
Jet Nestruev.
\newblock {\em Smooth manifolds and observables}, volume 220 of {\em Graduate Texts in Mathematics}.
\newblock Springer-Verlag, New York, russian edition, 2003.
\newblock Joint work of A. M.\ Astashov, A. B.\ Bocharov, S. V.\ Duzhin, A. B.\ Sossinsky, A. M.\ Vinogradov and M. M.\ Vinogradov.

\bibitem{BP}
Brent Pym.
\newblock {\em Poisson {S}tructures and {L}ie {A}lgebroids in {C}omplex {G}eometry}.
\newblock ProQuest LLC, Ann Arbor, MI, 2013.
\newblock Thesis (Ph.D.)--University of Toronto (Canada).

\bibitem{SR}
S.~Ramanan.
\newblock {\em Global calculus}, volume~65 of {\em Graduate Studies in Mathematics}.
\newblock American Mathematical Society, Providence, RI, 2005.

\bibitem{GR}
George~S. Rinehart.
\newblock Differential forms on general commutative algebras.
\newblock {\em Trans. Amer. Math. Soc.}, 108:195--222, 1963.

\bibitem{AS}
Abhishek Sarkar.
\newblock Cohomology of {L}ie algebroid over algebraic spaces.
\newblock {\em arXiv:2111.01735v8 [math.DG]}, 2024.

\bibitem{TS}
Travis Schedler.
\newblock Deformations of algebras in noncommutative geometry.
\newblock In {\em Noncommutative algebraic geometry}, volume~64 of {\em Math. Sci. Res. Inst. Publ.}, pages 71--165. Cambridge Univ. Press, New York, 2016.

\bibitem{Simpson}
Carlos~T. Simpson.
\newblock Moduli of representations of the fundamental group of a smooth projective variety. {I}.
\newblock {\em Inst. Hautes \'Etudes Sci. Publ. Math.}, (79):47--129, 1994.

\bibitem{Sridh}
R.~Sridharan.
\newblock Filtered algebras and representations of {L}ie algebras.
\newblock {\em Trans. Amer. Math. Soc.}, 100:530--550, 1961.

\bibitem{PT}
Pietro Tortella.
\newblock {$\Lambda$}-modules and holomorphic {L}ie algebroid connections.
\newblock {\em Cent. Eur. J. Math.}, 10(4):1422--1441, 2012.

\bibitem{JV}
Joel Villatoro.
\newblock On sheaves of {L}ie-{R}inehart algebras.
\newblock {\em arXiv:2010.15463v2 [math.DG]}, 2021.

\bibitem{CW}
Charles~A. Weibel.
\newblock {\em The {$K$}-book}, volume 145 of {\em Graduate Studies in Mathematics}.
\newblock American Mathematical Society, Providence, RI, 2013.
\newblock An introduction to algebraic $K$-theory.

\end{thebibliography}

\vspace{.1 cm}
{\bf Satyendra Kumar Mishra}, 
Department of Mathematics and Statistics,
Indian Institute of Technology (BHU) Varanasi, 
Varanasi 221005,
Uttar Pradesh, India.
e-mail: satyendrakr.mat@itbhu.ac.in\\

{\bf Abhishek Sarkar},
Department of Mathematics,
Indian Institute of Science Education and Research Pune, 
Pune 411008,
Maharashtra, India.
e-mail: abhisheksarkar49@gmail.com

	\end{document}